\documentclass[12pt,leqno]{article}

\usepackage[]{amsmath, amsthm, amsfonts}
\usepackage{amssymb}            
\usepackage[mathcal]{euscript}

\usepackage[matrix,arrow,curve,frame]{xy}       
\xymatrixcolsep{1.9pc}                          
\xymatrixrowsep{1.9pc}
\newdir{ >}{{}*!/-5pt/\dir{>}}

\numberwithin{equation}{section}

\newtheorem{theorem}{Theorem}[section]
\newtheorem{proposition}[theorem]{Proposition}

\newtheorem{corollary}[theorem]{Corollary}

\theoremstyle{definition}
\newtheorem{definition}[theorem]{Definition}
\newtheorem{pdef}[theorem]{Provisional Definition}
\newtheorem{example}[theorem]{Example}

\newtheorem{question}[theorem]{Question}

\theoremstyle{remark}
\newtheorem{remark}[theorem]{Remark}

\DeclareMathOperator{\Map}{Map}
\DeclareMathOperator{\Hom}{Hom}
\DeclareMathOperator{\Tot}{Tot}

\DeclareMathOperator*{\colim}{colim}

\newcommand{\bx}{\mathbin{\square}}

\renewcommand{\b}{\bullet}
\newcommand{\A}{{\cal A}}
\newcommand{\B}{{\cal B}}
\newcommand{\C}{{\cal C}}
\newcommand{\D}{{\cal D}}
\newcommand{\E}{{\cal E}}
\newcommand{\F}{{\cal F}}
\newcommand{\G}{{\cal G}}
\newcommand{\I}{{\cal I}}

\newcommand{\M}{{\cal M}}

\renewcommand{\S}{{\cal S}}
\newcommand{\T}{{\cal T}}
\newcommand{\U}{{\cal U}}
\renewcommand{\O}{{\cal O}}

\newcommand{\ds}{\displaystyle}

\begin{document}

\title{Operads and cosimplicial objects: an introduction.}
\author{James E. McClure and Jeffrey H. Smith%
\thanks{Both authors were supported by NSF grants.
This paper is based in part on lectures given by the first author at the Newton
Institute.}
\\Department of Mathematics, Purdue University  \\
150 N. University Street \\
West Lafayette, IN  47907-2067
}
\date{}
\maketitle

\section{Introduction.}

This paper is an introduction to the 
series of papers \cite{MS1,MS3,MS4,MS5}, in which we 
develop a combinatorial theory of certain important operads and their 
actions.\footnote{The relation between these papers and \cite{MS2} is 
explained in Remarks \ref{rev2} and \ref{rev3}.}  The operads we consider are
$A_\infty$ operads, $E_\infty$ operads, the little $n$-cubes operad and the
framed little disks operad.
Sections \ref{s2}, \ref{s6}
and \ref{s9}, which can be read independently, are an introduction to the
theory of operads. 

The reader is also referred to the very interesting papers 
of Batanin (\cite{Bat02,Bat03}), which treat similar questions from a 
categorical point of view.

Here is an outline of the paper.

In Section \ref{s2} we motivate the concept of non-symmetric operad.  We give
the definition of
$A_\infty$ space and state the 
characterization (up to weak equivalence) of loop
spaces: a space is weakly equivalent to a loop space if and only if it is a 
grouplike $A_\infty$ space.

In Section \ref{s3} we introduce the total space construction Tot for
cosimplicial spaces and the related conormalization construction 
for cosimplicial abelian groups.

In Section \ref{s4} we address the question: when is Tot of a cosimplicial
space an $A_\infty$ space?  We obtain a useful sufficient condition for this to
happen.

In Section \ref{s5} we reformulate the main result of Section \ref{s4} in a way
which is convenient for generalization.

In Section \ref{s6} we motivate the concept of operad.  We give the definition
of $E_\infty$ space and state the characterization (up to weak equivalence) of
infinite loop spaces: a space is weakly equivalent to an infinite loop space if
and only if it is a grouplike $E_\infty$ space.

Section \ref{s7} contains motivation for the main result of Section \ref{s8}.

In Section \ref{s8} we give a sufficient condition for Tot of a cosimplicial
space to be an $E_\infty$ space.

In Section \ref{s9} we introduce the little $n$-cubes operad $\C_n$.  Operads 
weakly equivalent to $\C_n$ are called $E_n$ operads.  We give the
characterization (up to weak equivalence) of $n$-fold loop spaces: a space is
weakly equivalent to an $n$-fold loop space if and only if it has a grouplike
action of an $E_n$ operad.

In Section \ref{s10} we give a sufficient condition for Tot of a cosimplicial
space to be an $E_n$ space.

In Section \ref{s11} we describe some category theory which is used in the
proof of the main theorem of Section \ref{s10}.

In Section \ref{s12} we outline the proof of the main theorem of Section 
\ref{s10}.  As a byproduct we get a new, combinatorial, description (up to 
weak equivalence) of $\C_n$.

In Section \ref{s13} we describe applications of the main result in Section
\ref{s10} to a certain space of knots and to topological Hochschild
cohomology.

In Section \ref{s14} we develop a combinatorial description of the framed 
little disks operad. 

In Section \ref{s15} we observe that the theory of 
Sections \ref{s4}, \ref{s5}, \ref{s8},
\ref{s10}, \ref{s12} and \ref{s14} remains
valid with spaces replaced by chain complexes.  In particular this leads to
concrete and explicit chain models for $\C_n$ and for the framed little
disks operad.

In Secion \ref{s16} we give some applications of the theory developed in
Section \ref{s15}; in particular we discuss Deligne's Hochschild cohomology
conjecture.

\section{Loop spaces and the little intervals operad.}

\label{s2}

Historically, the first use of operads was to give a precise
meaning to the idea that loop spaces are monoids up to higher homotopy.  In
this section we recall how this works.

The first step is to reformulate the concept of monoid in a way that is
amenable to generalization.  

\begin{proposition}
\label{1}
A monoid structure on a set $S$ determines and is determined by a family of 
maps
\[
M(k):S^k\to S
\]
for $k\geq 0$
(where $S^k$ denotes the $k$-fold Cartesian product)
such that 

(a) $M(1)$ is the identity map, and

(b) the set $\{M(k)\}_{k\geq 0}$ is closed under multivariable
composition.
\end{proposition}

\begin{proof}
If $S$ is a monoid with multiplication $M:S^2\to S$ and unit $e:S^0\to
S$ we define $M(0)$ to be $e$, $M(1)$ to be the identity map, and $M(k)$ to be
the iterated multiplication for $k\geq 2$.  The monoid axioms show that the
set $\{M(k)\}_{k\geq 0}$ is closed under multivariable
composition.

Conversely, if $S$ is a set with maps $M(k)$ satisfying (a) and (b) then $S$ is
a monoid with multiplication $M(2)$ and unit $M(0)$.
\end{proof}

Next let $Z$ be a based space with basepoint denoted by $*$.  We consider the 
space $\Omega Z$ of based loops on $Z$.  For each $r\in(0,1)$ there is a
multiplication
\[
M_r:(\Omega Z)^2 \to \Omega Z
\]
which takes a pair of loops $(\alpha,\beta)$ to the loop $M_r(\alpha,\beta)$
which is $\alpha$ (suitably rescaled) on the interval $[0,r]$ and $\beta$
(suitably rescaled) on the interval $[r,1]$. We represent the loop
$M_r(\alpha,\beta)$ by the picture

\noindent
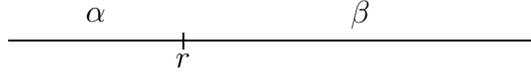
\mbox{}
\hfill
\begin{picture}(200,50)

\put(0,25){\line(1,0){200}}
\put(66,22){\line(0,1){6}}
\put(66,20){\makebox(0,0)[t]{$r$}}
\put(0,25){\makebox(66,20){$\alpha$}}
\put(66,25){\makebox(134,20){$\beta$}}

\end{picture}
\hfill
\mbox{}

%
%
%

\noindent
We write $*\in \Omega Z$ for the constant loop at the basepoint, 
which we represent by the picture

\noindent
\mbox{}
\hfill
\begin{picture}(200,50)

\put(0,25){\line(1,0){200}}
\put(0,25){\makebox(200,20){$*$}}

\end{picture}
\hfill
\mbox{}

\noindent
and we write
\[
e:(\Omega Z)^0 \to \Omega Z
\]
for the map whose image is $*$.

Motivated by Proposition \ref{1}, we consider the space $\M(k)$ of all maps 
\[
(\Omega Z)^k\to \Omega Z
\] 
that can be obtained by multivariable composition from the maps $M_r$ and $e$.
A typical example is the map in $\M(4)$ which takes a 4-tuple 
$(\alpha_1,\alpha_2,\alpha_3,\alpha_4)$ to the loop

%
%
%

\noindent
\mbox{}
\hfill
\begin{picture}(200,50)

\put(0,25){\line(1,0){200}}
\put(35,22){\line(0,1){6}}
\put(60,22){\line(0,1){6}}
\put(75,22){\line(0,1){6}}
\put(100,22){\line(0,1){6}}
\put(140,22){\line(0,1){6}}
\put(155,22){\line(0,1){6}}

\put(0,25){\makebox(35,15){$*$}}
\put(35,25){\makebox(25,15){$\alpha_1$}}
\put(60,25){\makebox(15,15){$\alpha_2$}}
\put(75,25){\makebox(25,15){$*$}}
\put(100,25){\makebox(40,15){$\alpha_3$}}
\put(140,25){\makebox(15,15){$*$}}
\put(155,25){\makebox(45,15){$\alpha_4$}}

\end{picture}
\hfill
\mbox{}

\noindent
In general, a map in $\M(k)$ is determined by $k$ closed intervals in $[0,1]$ 
with disjoint interiors (notice that, as in the example just given, the union 
of these intervals doesn't have to be all of $[0,1]$).  This motivates the
following definition:

\begin{definition}
\label{2}
Let $\A(k)$ be the set in which an element is a set of $k$ closed intervals in
$[0,1]$ with disjoint interiors (in particular, $\A(0)$ is a point).  Give 
$\A(k)$ the topology induced by the following imbedding of $\A(k)$ in ${\mathbb
R}^{2k}$: given $k$ closed intervals in $[0,1]$, list the $2k$ endpoints of 
the intervals in increasing order.
\end{definition}

What we have shown so far is that $\M(k)$ is homeomorphic to 
$\A(k)$.\footnote{In making this statement we must exclude the 
special case where the
path-component of the basepoint in $Z$ is a single point.  On the other hand,
the proposition which follows remains true in this case.}
Moreover, it is easy to see that each $\A(k)$ is contractible, so 
to sum up we have 

\begin{proposition}
\label{3}
If $Y=\Omega Z$ for some $Z$ then 
there is a family of subspaces $\M(k)\subset
\Map(Y^k,Y)$ 
such that

(a) $\M(1)$ contains the identity map,

(b) the family $\M=\{\M_k\}_k\geq 0$ is closed under multivariable composition,
and

(c) each $\M(k)$ is contractible.
\end{proposition}

The crucial fact about this situation is that Proposition \ref{3} has a
converse up to weak equivalence: if $Y$ is any connected space which has a 
family of contractible subspaces $\M(k)\subset
\Map(Y^k,Y)$ satisfying (a), (b) and (c) then $Y$ is weakly
equivalent to $\Omega Z$ for some $Z$ (this is a special case of Theorem 
\ref{10} below).  This gives us a way of recognizing that a space is a loop 
space (up to weak equivalence) without knowing in advance that it is a loop 
space.

Motivated by Propositions \ref{1} and \ref{3} we make a first attempt
at the definition of non-symmetric\footnote{The word ``non-symmetric'' refers to
the fact that we haven't yet used, or needed, the action of the symmetric 
group on $Y^k$; see Section \ref{s6}.} operad.

\begin{pdef}
\label{4}
A non-symmetric operad $\O$ is a collection of subspaces 
\[
\O(k)\subset \Map(Y^k,Y)\qquad k\geq 0
\]
(for some space $Y$) such that

(a) $\O(1)$ contains the identity map and

(b) the collection $\O$ is closed under multivariable composition.
\end{pdef}

\medskip

\noindent
{\bf Critique of Provisional Definition \ref{4}.}  This definition is formally
analogous to the nineteenth-century definition of a group as a family of
bijections of a set $S$, closed under composition and inverses.  The advantage
of such a definition is its concreteness and the ease with which our minds
assimilate it.  The disadvantage (in the case of groups) is that the set $S$ is
really external to the group.  The resolution of this difficulty (in the case 
of groups) was to split the original definition into two concepts: the 
concept of (abstract) group and the concept of group action.  

\medskip

Motivated by the Critique, we will split the Provisional Definition into two
concepts: the concept of (abstract) non-symmetric operad and the concept of
operad action. 

First observe that (in the situation of \ref{4}) the multivariable composition
operations in $\{\Map(Y^k,Y)\}_{k\geq 0}$ restrict to give maps
\[
\gamma:\O(k)\times\O(j_1)\times\cdots\times\O(j_k)\to \O(j_1+\cdots+j_k)
\]
for each choice of $k,j_1,\ldots,j_k\geq 0$.  The associativity property of
multivariable composition implies that the following diagram commutes for all
choices of $k,j_1,\ldots,j_k,\{i_{mn}\}_{m\leq k,\ n\leq j_m}$.
\begin{equation}
\label{5}
\xymatrix{
\O(k)\times\ds{\prod_{m=1}^k}
\Bigl(
\O(j_m)\times\prod_{n=1}^{j_m} \O(i_{mn})
\Bigr)
\ar[d]_{=}
\ar[r]^{1\times \gamma}
&
\O(k)\times \ds\prod_{m=1}^k \O(i_{m1}+\cdots+i_{mj_m})
\ar[dd]^{\gamma}
\\
\Bigl(
\O(k)\times \ds\prod_{m=1}^k \O(j_m)
\Bigr)
\times
\prod_{m,n}\O(i_{mn}) 
\ar[d]_-{\gamma\times 1}
&
\\
\O(j_1+\cdots+j_k)\times
\O(i_{11})\times\cdots\times
\O(i_{kj_k})
\ar[r]^-{\gamma}
&
\O(i_{11}+\cdots+i_{kj_k})
}
\end{equation}

\begin{definition}
\label{6}
A {\it non-symmetric operad} $\O$ is a collection of spaces $\{\O(k)\}_{k\geq
0}$ together with an element $1\in \O(1)$ and maps
\[
\gamma:\O(k)\times\O(j_1)\times\cdots\times\O(j_k)\to \O(j_1+\cdots+j_k)
\]
(for each choice of $k,j_1,\ldots,j_k\geq 0$)
such that 

(a) for each $k$ and each $s\in \O(k)$,
$\gamma(1,s)=s$ and $\gamma(s,1,\ldots,1)=s$, and

(b) Diagram \eqref{5} commutes for all choices of
$k,j_1,\ldots,j_k,i_{11},\ldots,i_{kj_k}$.
\end{definition}


\begin{remark}
\label{7}
Any collection $\O$ which satisfies Provisional Definition \ref{4} will
also satisfy Definition \ref{6}.  
\end{remark}

\begin{remark}
\label{7a} 
Here are some examples which will be important later.

(a) The collection $\A=\{\A(k)\}_{k\geq 0}$ defined in Definition \ref{2} is 
a non-symmetric operad (by Remark \ref{7}); it is called the {\it little 
intervals non-symmetric operad}.
It is instructive to work out 
the explicit description of the maps $\gamma$ in this case (cf.\ Section
\ref{s9}).

(b) If $Y$ is any space the collection 
\[
\{ \Map(Y^k,Y) \}_{k\geq 0},
\]
with its usual multivariable composition, is called the 
{\it endomorphism operad of\/} $Y$ and denoted $\E_Y$.  

(c) More generally, if $\U$ is any topological category with a monoidal product
$\bx$ (see \cite[Section VII.1]{MacLane})
and $D$ is an object of $\U$ then the collection of spaces
\[
\{ \Hom_\U(D^{\bx k},D) \}_{k\geq 0}
\]
with the evident multivariable composition is a non-symmetric 
operad.  The proof is left as an exercise for the reader; see Section \ref{s11}
for a hint.

(d) With $\U$ as above, note that $\bx$ is also a monoidal product for 
$\U^{\text{op}}$. 
Applying part (c) to $\U^{\text{op}}$ gives a
non-symmetric operad whose $k$-th space is
\[
\Hom_{\U^{\text{op}}}(D^{\bx k},D)=
\Hom_\U(D,D^{\bx k})
\] 
\end{remark}

Next we formulate the concept of operad action.  First observe that in the
situation of Provisional Definition \ref{4} the evaluation maps
\[
\Map(Y^k,Y)\times Y^k \to Y
\]
restrict to maps
\[
\theta: \O(k)\times Y^k\to Y
\]
and that the diagram 
\begin{equation}
\label{8}
\xymatrix{
\O(k)\times
\Bigl(
\ds\prod_{m=1}^k \O(j_m)
\Bigr)
\times Y^{j_1+\cdots+j_k}
\ar[r]^-{=}
\ar[dd]_-{\gamma\times 1}
&
\O(k)\times\ds\prod_{m=1}^k
\bigl(
\O(j_m)\times Y^{j_m}
\bigr)
\ar[d]^-{1\times\prod\theta}
\\
&
\O(k)\times Y^k
\ar[d]^-{\theta}
\\
\O(j_1+\cdots+j_k)\times Y^{j_1+\cdots+j_k}
\ar[r]^-{\theta}
&
Y
}
\end{equation}
commutes for all $k,j_1,\ldots,j_k\geq 0$.

\begin{definition}
\label{8a}
Let $\O$ be a non-symmetric operad and let $Y$ be a space.  An action of 
$\O$ on $Y$ consists of a map
\[
\theta: \O(k)\times Y^k\to Y
\]
for each $k\geq 0$ such that

(a) $\theta(1,x)=x$ for all $x\in Y$, and

(b) Diagram \eqref{8} commutes for all $k,j_1,\ldots,j_k\geq 0$.
\end{definition}

\begin{example}
\label{9}
(a) The non-symmetric operad $\A$ mentioned in Remark \ref{7a}(a) acts on 
$\Omega Z$ for any space $Z$.

(b) The endomorphism operad $\E_Y$ acts on $Y$.
\end{example}

We conclude this section by stating the most general converse to Proposition 
\ref{3}.  First we need two definitions.

\begin{definition}
\label{9a}
An $A_\infty$ operad is a non-symmetric operad $\O$ for which each space
$\O(k)$ is weakly equivalent to a point.
\end{definition}

For example, the non-symmetric operad $\A$ of Remark \ref{7a}(a) is an
$A_\infty$ operad.

Notice the relationship between this definition and Proposition \ref{1}: the
one-point sets $\{M(k)\}$ in Proposition \ref{1} are replaced by contractible
spaces in Definition \ref{9a}.

Now let $Y$ be a space with an action 
of an $A_\infty$ operad $\O$.  Because $\O(2)$ and $\O(0)$ are connected, the 
maps
\[
\theta:\O(2)\times Y^2\to Y
\]
and 
\[
\theta:\O(0)\times Y^0\to Y
\]
induce a monoid structure on $\pi_0 Y$.  

\begin{definition}
The action of $\O$ on $Y$ is {\it grouplike} if the monoid $\pi_0 Y$ is a 
group.
\end{definition}

For example, the action in Example \ref{9}(a) is grouplike.  Also, if $Y$ is 
connected then all actions are grouplike.

\begin{theorem}
\label{10}
$Y$ is weakly 
equivalent to $\Omega Z$ for some space $Z$ $\Longleftrightarrow$
$Y$ has a grouplike action of an $A_\infty$ operad. 
\end{theorem}

\begin{remark}
This theorem developed gradually during the period from 1960 to 1974.  
In the $\Longleftarrow$ direction,
the first version was proved by Stasheff \cite{Stasheff}, assuming that $Y$ is 
connected and using a particular non-symmetric operad, now called the 
Stasheff operad (but the concept of operad hadn't yet been defined at that 
time).  Boardman and Vogt proved the $\Longleftarrow$ direction for general 
$A_\infty$ operads (except that they used PROP's instead of operads), but 
still assuming $Y$ connected, in \cite{BV1,BV2}. 
May defined the concept of 
operad in \cite{MayG} and proved the $\Longleftarrow$ direction for connected 
$Y$; he proved the general version (for group-complete actions) in \cite{MayE}.
The $\Longrightarrow$ direction (for PROP's, which implies the result for
operads) is due to Boardman and Vogt \cite{BV1,BV2}.
\end{remark}

\section{Cosimplicial objects and totalization.}
\label{s3}

Theorem \ref{10} leads to the question of how we can tell when a space $Y$ has
an action of an $A_\infty$ operad.  In the next section we will give an 
answer to this question in the important special case where $Y$ is the total 
space of a cosimplicial space $X^\b$.  In this section we pause for some
background about cosimplicial objects.

Throughout this paper we will use the following conventions for cosimplicial
objects.

\begin{definition}
\label{11}
(a) Define $\Delta$ to be the category of nonempty finite totally ordered
sets (this is equivalent to the category usually called $\Delta$).  Define
$[m]$ to be the finite totally ordered set $\{0,\ldots,m\}$. Define
\[
d^i:[m]\to [m+1] \qquad 0\leq i \leq m+1
\]
to be the unique ordered injection whose image does not contain $i$,
and define
\[
s^i:[m]\to [m-1] \qquad 0\leq i \leq m-1
\]
to be the unique ordered surjection for which the inverse image of $i$ contains
two points.

(b) Given a category $\C$, 
a cosimplicial object $X^\b$ in $C$ is a functor from $\Delta$ to $\C$.
If $S$ is a nonempty finite totally ordered set then 
$X^S$ will denote the value of $X^\b$ at $S$, except that we write $X^m$ 
instead of $X^{[m]}$.  The maps
\[
X^m\to X^{m+1}
\]
induced by the $d^i$ are called {\it coface maps} and the maps
\[
X^m\to X^{m-1}
\]
induced by the $s^i$ are called {\it codegeneracy maps}.
\end{definition}

Note that every object in $\Delta$ has a unique isomorphism to an object of the
form $[m]$, so we can specify a cosimplicial object by giving its value on the
objects $[m]$ (together with the coface and codegeneracy maps).  For example:

\begin{definition}
$\Delta^\b$ is the cosimplicial space whose value at $[m]$ is the simplex 
$\Delta^m$, with the usual coface and codegeneracy maps.
\end{definition}

Next we define the cosimplicial analog of geometric realization. 
First recall (for example, from \cite[Example 2.4(3)]{HV}) that the 
geometric realization of a simplicial space $U_\b$ is a
tensor product over $\Delta$ (also called a coend):
\[
|U_\b|=U_\b \otimes_{\Delta} \Delta^\b
\] 
When we change the variance from simplicial to cosimplicial it is natural to
replace $\otimes_\Delta$ by $\Hom_\Delta$, which leads us to the following 
definition.

\begin{definition}
\label{11aaaa}
Let $X^\b$ be a cosimplicial space.  The {\it total space} of $X^\b$, denoted
$\Tot(X^\b)$, is the space of cosimplicial maps $\Hom_\Delta(\Delta^\b,X^\b)$.
\end{definition}

Here's a more explicit description:
a point in $\Tot(X^\bullet)$ is a sequence
\[
\alpha_0:\Delta^0\to X^0, \alpha_1:\Delta^1\rightarrow X^1, 
\alpha_2:\Delta^2\rightarrow X^2, \ldots
\]
which is {\it consistent}, i.e.,
\[
d^i \circ \alpha_n = \alpha_{n+1}\circ d^i
\]
and
\[
s^i \circ \alpha_n = \alpha_{n-1} \circ s^i
\]
for all $i$.

\begin{example}
\label{11aaa}
Given a based space $Z$ with basepoint $*$, 
we define a cosimplicial space $F^\b Z$ whose total space is $\Omega Z$ ($F^\b
Z$ is called the geometric cobar construction on $Z$). The $m$-th space
$F^m Z$ is the Cartesian product $Z^m$.  The coface 
\[
d^i:F^m Z\to F^{m+1} Z
\]
is defined by
\[
d^i(z_1,\ldots,z_m)=\left\{
\begin{array}{ll}
(*,z_1,\ldots,z_m) & \mbox{if $i=0$} \\
(z_1,\ldots, z_i,z_i,\ldots,z_m) & \mbox{if $1\leq i\leq m$} \\
(z_1,\ldots,z_m,*) & \mbox{if $i=m+1$}.
\end{array}
\right.
\]
and the codegeneracy $s^i:F^m Z\to F^{m-1} Z$ deletes the $(i-1)$-st
coordinate.  The proof that
$\Tot(F^\b Z)$ is homeomorphic to $\Omega Z$ is left as an
exercise for the reader. (Hint: if $m>1$ then the map
\[
\prod_{i=0}^{m-1} s_i :F^m Z\to \prod_{i=0}^{m-1} F^{m-1} Z
\]
is a monomorphism).
\end{example}

We will also consider cosimplicial abelian groups. 

\begin{example}
\label{11aa}
Let $W$ be a space and define $S^\b W$ to be the cosimplicial abelian group 
$\Map(S_\b W, {\mathbb Z})$, where $S_\b W$ is the usual simplicial set 
associated to $W$ (\cite[Example 1.28]{Curtis}) and Map means maps of sets.
\end{example}

Next we define the analog of Tot in this context. 
Let $\Delta_{\text{simp}}^m$ be the standard simplicial model of $\Delta^m$
(\cite[Example 1.4]{Curtis}). 

\begin{definition}
\label{S1}
Let $\Delta^\b_*$ denote the cosimplicial chain 
complex which in degree $m$ is the normalized chain complex 
(\cite[pages 265--266]{Weibel})
of $\Delta_{\text{simp}}^m$.  
\end{definition}

\begin{definition}
\label{11a}
Let $A^\b$ be a cosimplicial abelian group.  The {\it conormalization}%
\footnote{This functor is usually called normalization, but it seems
desirable to have separate names for the cosimplicial and simplicial 
versions of normalization,
analogous to the usual distinction between Tot and geometric realization.}
of $A^\b$, denoted ${\mathbf C}(A^\b)$, is the cochain complex
\[
\Hom_\Delta (\Delta^\b_*,A^\b) \subset
\prod_{m=0}^\infty \Hom(\Delta^m_*,A^m).
\]
Here $\Hom_\Delta$ is Hom in the category of cosimplicial graded abelian 
groups (with $A^m$ concentrated in dimension 0), and the differential
is induced by the differentials of the $\Delta^m_*$.
\end{definition}

\begin{remark}
\label{12}
Here are two concrete descriptions of ${\mathbf C}(A^\b)$; they are dual to 
the two
standard ways of describing the normalization of a simplicial abelian group
(\cite[pages 265--266]{Weibel}).  The proof that they agree with Definition
\ref{11a} is left to the reader.

(a) 
Let ${\mathbf C}'(A^\b)$ be the cochain complex whose $m$-th group is
the intersection of the kernels of the codegeneracies $s^i:A^m\to A^{m-1}$
and whose differential is $\sum (-1)^i d^i$.  Then ${\mathbf C}(A^\b)$ is 
isomorphic to ${\mathbf C}'(A^\b)$.

(b)
Let ${\mathbf C}''(A^\b)$ be the cochain complex whose $m$-th group is the 
cokernel of
\[
\bigoplus_{i>0} \,d^i: \bigoplus_{i>0}\, A^{m-1}\to A^m
\]
and whose differential is induced by $d^0$.  Then ${\mathbf C}(A^\b)$ is 
isomorphic to ${\mathbf C}''(A^\b)$
\end{remark}

\begin{example}
\label{12a}
The conormalization of $S^\b W$ (Example \ref{11aa}) is the complex 
of singular cochains that vanish on all degenerate singular chains.
This is what is usually called the normalized
singular cochain complex of $W$; we will denote it by 
$s^* W$.  
\end{example}

\section{A sufficient condition for $\Tot(X^\b)$ to be an $A_\infty$ space.}
\label{s4}

\begin{definition}
An $A_\infty$ space is a space with an action of an $A_\infty$ operad.
\end{definition}

Let $Z$ be a based space and let $F^\b Z$ be the cosimplicial space defined in 
Example \ref{11aaa}.  Then $\Tot(F^\b Z)$ is homeomorphic to $\Omega Z$ and in
particular (as we have seen in Section \ref{s2}) it is an $A_\infty$ space.
This leads us to the question:

\begin{question}
\label{13}
For what other cosimplicial spaces is Tot an $A_\infty$ space?
\end{question}

As we have seen in Section \ref{s3}, Tot is analogous to conormalization, so we
can gain insight into Question \ref{13} by examining a cosimplicial abelian
group whose conormalization has a multiplicative structure,
namely $S^\b W$ (see Example \ref{11aa}). The conormalization of $S^\b W$ is 
$s^* W$ (see Example \ref{12a}), and $s^* W$ has an associative multiplication,
the cup product, given by the usual Alexander-Whitney formula
\[
(x\smallsmile y)(\sigma)=x(\sigma(0,\ldots,p))\cdot
y(\sigma(p,\ldots,p+q));
\]
here $x$ has degree $p$, $y$ has degree $q$, $\sigma$ is in
$S_{p+q}W$, $\cdot$ is multiplication in $\mathbb Z$, and
$\sigma(0,\ldots,p)$ (resp.,
$\sigma(p,\ldots,p+q)$)
is the restriction of $\sigma$ to the subsimplex of
$\Delta^{p+q}$ spanned by the vertices $0,\ldots,p$ (resp., $p,\ldots,p+q$).
The key point for our purpose is that the same formula defines a map
\begin{equation}
\label{rev1}
\smallsmile:S^p W\times S^q W \to S^{p+q} W
\end{equation}
and we can examine the relation between
$\smallsmile$ and the coface and codegeneracy maps 
of $S^\b W$. This relation is
given by the following formulas:
\begin{equation}
\label{14}
d^i (x\smallsmile y)=\left\{
\begin{array}{ll}
d^i x \smallsmile y & \mbox{if $i\leq p$} \\
x\smallsmile d^{i-p}y & \mbox{if $i>p$}
\end{array}
\right.
\end{equation}
\begin{equation}
\label{15}
d^{p+1}x\smallsmile y=x\smallsmile d^0 y
\end{equation}
\begin{equation}
\label{16}
s^i(x\smallsmile y)=\left\{
\begin{array}{ll}
s^i x\smallsmile y & \mbox{if $i\leq p-1$} \\
x \smallsmile s^{i-p}y & \mbox{if $i \geq p$}
\end{array}
\right.
\end{equation}

Next we observe that the cosimplicial space
$F^\b Z$ has the same kind of structure as $S^\b W$: if we define
\[
\smallsmile: F^p Z\times F^ q Z \to F^{p+q} Z
\]
to be the obvious juxtaposition map 
\[
Z^p\times Z^q \to
Z^{p+q}
\]
then $\smallsmile$ satisfies \eqref{14}, 
\eqref{15} and \eqref{16}.  Moreover, it is associative:
\[
(x\smallsmile y)\smallsmile z=z\smallsmile (y\smallsmile z),
\]
and unital:
there is an element $e\in F^0 Z$ (namely the basepoint) such that
\[
x\smallsmile e=e\smallsmile x=x
\]
for all $x$.
This suggests that, as a way of answering Question \ref{13}, we consider
cosimplicial spaces having the same kind of structure as $S^\b W$ or $F^\b Z$:

\begin{theorem}
\label{17}
If $X^\b$ is a cosimplicial space with a cup product
\[
\smallsmile:X^p\times X^q \to X^{p+q}
\]
which is associative and unital and satisfies
\eqref{14}, \eqref{15} and \eqref{16} then $\Tot(X^\b)$ is an $A_\infty$ 
space.
\end{theorem}

\begin{remark}
(a) This result is due to Batanin \cite[Theorems 5.1 and 5.2]{Bat98} with a 
simplified proof by us \cite[Section 3]{MS3}.

(b) Theorem \ref{17} gives a sufficient but not a necessary condition for
$\Tot(X^\b)$ to be an $A_\infty$ space.  However, we expect that any 
$A_\infty$ space is weakly equivalent to one produced by Theorem \ref{17} (in 
fact it is likely that $\Tot$ induces a Quillen equivalence between cosimplicial
spaces satisfying the hypothesis of Theorem \ref{17} and $A_\infty$ spaces).
\end{remark}

The remainder of this section gives an outline of the proof of Theorem 
\ref{17}; for details see \cite[Section 3]{MS3}.

The first step in the proof is:

\begin{proposition}
The category of cosimplicial spaces has a monoidal structure $\bx$ with the
property that $X^\b$ satisfies the hypothesis of Theorem \ref{17} if and 
only if it is a $\bx$-monoid.
\end{proposition}

This is due to Batanin \cite{Bat93}.  The definition of $\bx$ is modeled on
equations \eqref{14}, \eqref{15} and \eqref{16}:
$X^\b\bx Y^\b$ is
the cosimplicial space whose $m$-th space is
\[
\left(\coprod_{p+q=m} X^p\times Y^q\right)/\sim
\]
where $\sim$ is the equivalence relation generated by $(x,d^0 y)\sim
(d^{|x|+1}x,y)$.
The coface maps are defined by
\[
d^i(x,y)=\left\{
\begin{array}{l}
(d^i x,y)\mbox{ if $i\leq |x|$} \\
(x,d^{i-|x|}y)\mbox{ if $i>|x|$}
\end{array}
\right.
\]
and the codegeneracy maps by
\[
s^i(x,y)=\left\{
\begin{array}{l}
(s^i x,y)\mbox{ if $i\leq |x|-1$} \\
(x,s^{i-|x|}y)\mbox{ if $i\geq |x|$}
\end{array}
\right.
\]

Next we apply Remark \ref{7a}(d) with $D=\Delta^\b$ to get a non-symmetric 
operad $\B$.  The space $\B(k)$ is
\[
\Hom_\Delta(\Delta^\b,(\Delta^\b)^{\bx k})
\]
The composition maps $\gamma$ are defined as follows:
if $f\in \B(k)$ and $g_i\in \B(j_i)$ for $1\leq i \leq k$ then
$\gamma(f,g_1,\ldots,g_k)\in \B(j_1+\cdots+j_k)$
is the composite
\[
\Delta^\b\xrightarrow{f}(\Delta^\b)^{\bx k}
\xrightarrow{g_1\bx\cdots\bx g_k}
(\Delta^\b)^{\bx(j_1+\cdots+j_k)}
\]

Now let $X^\b$ be a $\bx$-monoid.  We define an action of $\B$ on $\Tot(X^\b)$
by letting
\[
\theta:\B(k)\times (\Tot(X^\b))^k \to \Tot(X^\b)
\]
take $(f,\tau_1,\ldots,\tau_k)$ to the composite
\[
\Delta^\b\xrightarrow{f}(\Delta^\b)^{\bx k}
\xrightarrow{\tau_1\bx\cdots\bx \tau_k} (X^\b)^{\bx k}
\xrightarrow{\mu} X^\b
\]
where $\mu$ is the monoidal structure map of $X^\b$.

To complete the proof of Theorem \ref{17} it only remains to show that each
$\B(k)$ is contractible.  This is an easy consequence of the fact (due to
Grayson \cite{Grayson}) that $(\Delta^\b)^{\bx k}$ is isomorphic 
as a cosimplicial space to $\Delta^\b$.  See \cite[Section 3]{MS3} for details.

\section{A reformulation.}
\label{s5}

Our next goal is to generalize Theorem \ref{17}.  However, it turns out that 
the analogs of equations \eqref{14}, \eqref{15} and \eqref{16} for the
situations we will be considering are rather complicated and inconvenient, so 
we pause to reformulate Theorem \ref{17} in a way that is more amenable to
generalization.

Let us return to the motivating example $S^\b W$. Define 
\[
\sqcup: S^p W \otimes S^q W \to S^{p+q+1} W
\]
by 
\[
(x\sqcup y)(\sigma)=x(\sigma(0,\ldots,p))\cdot
y(\sigma(p+1,\ldots,p+q+1))
\]
for $\sigma\in S_{p+q+1} W$.
Note that, in contrast to the cup product, the vertex $p$ is not repeated in
the formula for $\sqcup$.

This operation is related to the coface and codegeneracy operations in 
$S^\b W$
by the following equations:
\begin{equation}
\label{18}
d^i (x\sqcup y) =
\left\{
\begin{array}{ll}
d^i x \sqcup y & \mbox{if $i\leq p+1$} \\
x\sqcup d^{i-p-2}y & \mbox{if $i>p+1$}
\end{array}
\right.
\end{equation}
\begin{equation}
\label{19}
s^i (x\sqcup y) =
\left\{
\begin{array}{ll}
s^i x \sqcup y & \mbox{if $i< p$} \\
x\sqcup s^{i-p-1}y & \mbox{if $i>p$}
\end{array}
\right.
\end{equation}
Note that there is no analog for $\sqcup$ of equation \eqref{15}.
$\sqcup$ is associative:
\[
x\sqcup (y\sqcup z)=
(x\sqcup y)\sqcup z
\]
and unital: there exists $e\in X^0$ with
\begin{equation}
\label{21}
s^p(x\sqcup e)=s^0(e\sqcup x)=x.
\end{equation}

The operations $\smallsmile$ and $\sqcup$ determine each other:
\begin{equation}
\label{22}
x\sqcup y = (d^{p+1}x)\smallsmile y =x\smallsmile d^0 y
\end{equation}
\begin{equation}
\label{23}
x \smallsmile y =s^p(x \sqcup y)
\end{equation}

Now let $X^\b$ be a cosimplicial space.  If $X^\b$ has an operation
\[
\sqcup: X^{p}\times X^q \to X^{p+q+1}
\]
which is associative and unital (in the sense of equation \eqref{21})
and satisfies \eqref{18} and \eqref{19} then the operation $\smallsmile$
defined by equation \eqref{23} satisfies the hypothesis of Theorem \ref{17} 
(the verification is left to the reader).  Conversely, if $X^\b$ has a cup
product satisfying the hypothesis of Theorem \ref{17} then the $\sqcup$
product defined by \eqref{22} is associative, unital and satisfies 
\eqref{18} and \eqref{19}.  To sum up:

\begin{proposition}
$X^\b$ satisfies the hypothesis of Theorem \ref{17} if and only it has a
product 
\[
\sqcup:X^p\times X^q \to X^{p+q+1}
\]
which is associative and unital and satisfies \eqref{18} and \eqref{19}.
\end{proposition}

\begin{corollary}
\label{23a}
If $X^\b$ is a cosimplicial space with a product
\[
\sqcup:X^p\times X^q \to X^{p+q+1}
\]
which is associative and unital 
and satisfies
\eqref{18} and \eqref{19}
then $\Tot(X^\b)$ is an $A_\infty$
space.
\end{corollary}

\section{Operads.}
\label{s6}

As we have seen in Section \ref{s2}, the ``associativity up to higher 
homotopy'' of the multiplication on $\Omega Z$ can be formulated rigorously 
as the action of an $A_\infty$ operad on $\Omega Z$.  We would now like to 
give an analogous formulation of ``commutativity up to higher homotopy.''   A 
multiplication is commutative in the ordinary sense if it is invariant under 
permutations of the factors; this suggests that we add symmetric-group 
actions to Definition \ref{6}.  We begin with a provisional form
of the definition.

\begin{pdef}
\label{24}
An operad $\O$ is a collection of subspaces
\[
\O(k)\subset \Map(Y^k,Y)\qquad k\geq 0
\]
(for some space $Y$) such that

(a) $\O(1)$ contains the identity map, 

(b) the collection $\O$ is closed under multivariable composition, and

(c) each $\O(k)$ is closed under the permutation action of the symmetric group
$\Sigma_k$.
\end{pdef}

As in Section \ref{s2} we split the provisional definition into 
the concept of (abstract) operad and the concept of operad action.

To formulate the definition of operad, we need to know the relation between the
action of $\Sigma_k$ and the multivariable composition maps $\gamma$ in
Provisional Definition \ref{24}.  This is left as an exercise for the reader;
the answer in given in \cite[Definition 1.1(c)]{MayG}.

\begin{definition}
\label{24a}
An operad is a non-symmetric operad $\O$ together with, for each $k$, a right
action of $\Sigma_k$ satisfying the formulas of \cite[Definition
1.1(c)]{MayG}.
\end{definition}

\begin{remark}
\label{25}
(a) Let $Y$ be a space.  The endomorphism operad $\E_Y$ (Remark \ref{7a}(b))
is an operad, where the $k$-th space $\Map(Y^k,Y)$ is given the obvious right
action of $\Sigma_k$.

(b) If $\U$ is a topological category with a {\it symmetric} monoidal product
$\bx$ (see \cite[Section VII.7]{MacLane}) then the non-symmetric 
operads in Remarks \ref{7a}(c) and (d) are operads, with the obvious 
$\Sigma_k$ actions on $\Hom_\U(D^{\bx k},D)$ and $\Hom_\U(D,D^{\bx k})$.

(c) If $\O$ is a non-symmetric operad we can define an operad $\O'$ by 
\[
\O'(k)=\O(k)\times \Sigma_k
\]
with the obvious right $\Sigma_k$ action; the definition of the composition 
maps $\gamma$ for $\O'$ is left as an exercise.  We call $\O'$ the operad 
generated by $\O$.
\end{remark}

\begin{definition}
Let $\O$ be an operad and let $Y$ be a space.  An action of $\O$ on $Y$ is an
action of the underlying non-symmetric operad with the property that each map 
\[
\theta: \O(k)\times Y^k\to Y
\]
factors through $\O(k)\times_{\Sigma_k} Y^k$.
\end{definition}

\begin{remark}
If $Y$ is a space then $\E_Y$ acts on $Y$.
\end{remark}

\begin{definition}
An $E_\infty$ operad is an operad $\O$ for which each space $\O(k)$ is weakly 
equivalent to a point.\footnote{For technical reasons, it is usual to require
in addition that the action of each $\Sigma_k$ should be free. The operad 
$\D$ that we construct in Section \ref{s8} does have this property.}
\end{definition}

A space with an action of an $E_\infty$ operad should be thought of as
``commutative up to all higher homotopies.''

The analog of Theorem \ref{10} in this setting is a statement about ``infinite
loop spaces.''  Recall that an infinite loop space is a space $X$ for which
there exists a sequence $X_1, X_2, \ldots$ with $X$ homeomorphic to $\Omega
X_1$ and $X_i$ homeomorphic to $\Omega X_{i+1}$ for all $i$ (thus an infinite
loop space is the zeroth space of a spectrum).

\begin{theorem}
\label{26}
$Y$ is weakly
equivalent to an infinite loop space
$\Longleftrightarrow$
$Y$ has a grouplike action of an $E_\infty$ operad. 
\end{theorem}

The $\Longrightarrow$ direction, and the 
$\Longleftarrow$ direction for connected $Y$, 
are due to Boardman and Vogt 
\cite{BV1,BV2}.  May gave a simpler proof of the $\Longleftarrow$ direction for
connected $Y$ in \cite{MayG} and proved the general case in
\cite{MayE}.

\section{A family of cochain operations.}

\label{s7}

We want to give an $E_\infty$ analog of Corollary \ref{23a}. 
In this section we prepare the way by returning to the motivating example, 
$S^\b W$, and defining a family of operations that generalizes $\sqcup$.  The 
idea is that the definition of $\sqcup$ is based on the partition of the 
set $\{0,\ldots,p+q+1\}$ into $\{0,\ldots,p\}$ and $\{p+1,\ldots,p+q+1\}$; we 
can produce more operations by using other partitions.

First we need some notation.
Recall that we have defined $\Delta$ to be the category of nonempty 
finite totally ordered sets $T$.  
For $T\in \Delta$ we define $\Delta^T$ to be the convex hull of $T$ (in
particular, $\Delta^{[m]}$ is the usual $\Delta^m$).  We define $S_T W$ to be
the set of all continuous maps $\Delta^T\to W$ (in particular, $S_{[m]} W$ is
what we have been calling $S_m W$) and $S^T W$ to be 
$\Map(S_T W, {\mathbb Z})$ (so $S^{[m]} W$ is the same as $S^m W$).

\begin{definition}
Given a map $\sigma:\Delta^m\to W$ and a subset $U$ of $T$, let 
$\sigma(U)$ be the restriction of $\sigma$ to $\Delta^U$.
\end{definition}

Now observe that a partition of $T$
into two pieces is the same thing as a surjective function
$f:T \to \{1,2\}$.

\begin{definition}
\label{26a}
Given a surjection $f:T\to\{1,2\}$, define a natural transformation
\[
\langle f \rangle: S^{f^{-1}(1)} W\times S^{f^{-1}(2)} W
\to S^T W
\]
by the equation
\[
\langle f \rangle(x,y)(\sigma)=
x(\sigma(f^{-1}(1)))\cdot y(\sigma(f^{-1}(2)))
\]
for $\sigma\in S_{T} W$; here $\cdot$ is multiplication in $\mathbb Z$.
\end{definition}

\begin{remark}
\label{27}
If $f$ is the function $\{0,\ldots,p+q+1\}\to \{1,2\}$ that takes 
$\{0,\ldots,p\}$ to 1 and $\{p+1,\ldots,p+q+1\}$
to 2 then $\langle f \rangle$ is $\sqcup$.
\end{remark}

Next we describe the relation between the operations $\langle f\rangle$ and 
the cosimplicial structure maps of $S^\b W$. 

\begin{proposition} 
\label{28}
Let 
\[
\xymatrix{
T
\ar[rr]^{\phi}
\ar[dr]_-f
&&
T'
\ar[dl]^-g
\\
&
\{1,2\}
&
}
\]
be a commutative diagram, where 
$\phi$ is a map in $\Delta$ (i.e., an order-preserving map).  For $i=1,2$ let 
\[
\phi_i:f^{-1}(i)\to g^{-1}(i)
\]
be the restriction of $\phi$.  

Then the diagram
\[
\xymatrix{
S^{f^{-1}(1)}W\times S^{f^{-1}(2)}W
\ar[d]_{(\phi_1)_*\times (\phi_2)_*}
\ar[r]^-{\langle f\rangle}
&
S^T W
\ar[d]^{\phi_*}
\\
S^{g^{-1}(1)}W\times S^{g^{-1}(2)}W
\ar[r]^-{\langle g\rangle}
&
S^{T'}W
}
\]
commutes.  
\end{proposition}

The proof is an immediate consequence of the definitions.  In the special 
case of Remark \ref{27} we recover equations \eqref{18} and \eqref{19}.

Next we formulate the commutativity, associativity and unitality properties of
the $\langle f\rangle$ operations.  Commutativity is easy:

\begin{proposition}
\label{29}
The diagram
\[
\xymatrix{
S^{f^{-1}(1)} W\times S^{f^{-1}(2)} W
\ar[d]_{\tau}
\ar[r]^-{\langle f\rangle}
&
S^T W
\ar[d]^{=}
\\
S^{f^{-1}(2)} W\times S^{f^{-1}(1)} W
\ar[r]^-{\langle t\circ f\rangle}
&
S^T W
}
\]
commutes, where $\tau$ is the switch map and $t$ is the transposition of
$\{1,2\}$.
\end{proposition}

For the associativity condition we need some notation.  
Define
\[
\alpha:\{1,2,3\}\to\{1,2\} 
\]
by $\alpha(1)=1,\alpha(2)=1,\alpha(3)=2$ and 
\[
\beta:\{1,2,3\}\to\{1,2\}
\]
by $\beta(1)=1,\beta(2)=2,\beta(3)=2$.
Given a surjection $g:T\to \{1,2,3\}$
let $g_1$ be the restriction of $g$ to $g^{-1}\{1,2\}$ and let $g_2$ be the
restriction of $g$ to $g^{-1}\{2,3\}$.

\begin{proposition}
\label{30}
With the notation above, the diagram 
\[
\xymatrix{
S^{g^{-1}(1)} W \times S^{g^{-1}(2)} W \times S^{g^{-1}(3)} W
\ar[d]_{1\times \langle g_2\rangle}
\ar[r]^-{\langle g_1\rangle\times 1}
&
S^{g^{-1}\{1,2\}} W \times S^{g^{-1}(3)} W
\ar[d]^{\langle \alpha \circ g\rangle}
\\
S^{g^{-1}(1)} W \times S^{g^{-1}\{2,3\}} W
\ar[r]^-{\langle \beta\circ g\rangle}
&
S^T W
}
\]
commutes for every choice of $T$ and $g$.
\end{proposition}

Again, the proof is immediate from the definitions.

For the unital property we need to extend $S^\b W$ to the category 
of all finite totally ordered sets, including the empty set.  

\begin{definition}
(a) Define $\Delta_+$ to be the category of finite totally ordered sets.

(b) Given a category $\C$,
an {\it augmented cosimplicial object} in $C$ is a functor from $\Delta_+$ to 
$\C$.

(c) Extend $\Delta^\b$ to a functor on $\Delta_+$ by defining
$\Delta^\emptyset=\emptyset$.

(d) Define $S_\b W$ as a functor from $\Delta_+^{\text{op}}$ to sets by
$S_T W=\Map(\Delta^T, W)$; in particular $S_\emptyset W$ is a point.

(e) Define $S^\b W$ as a functor from $\Delta_+$ to abelian groups by
$S^T W=\Map(S_T W,{\mathbb Z})$; in particular $S^\emptyset W$ is isomorphic to
$\mathbb Z$.
\end{definition}

With these conventions, Definition \ref{26a} makes sense when $f$ is not 
surjective, and Propositions \ref{28} (with $\Delta$ replaced by $\Delta_+$), 
\ref{29} and \ref{30} are still valid in this slightly more general context.

Now let $\varepsilon\in S^\emptyset W$ be the element corresponding to 
$1\in {\mathbb Z}$.

\begin{proposition}
\label{31}
If $f:T\to\{1,2\}$ takes all of $T$ to 1 then $\langle f\rangle 
(x,\varepsilon)=x$ for all $x$ and if $f$ takes all of $T$ to 2 then $\langle 
f\rangle (\varepsilon,x)=x$ for all $x$.
\end{proposition}

\section{A sufficient condition for $\Tot(X^\b)$ to be an $E_\infty$ space.}

\label{s8}

\begin{definition}
An $E_\infty$ space is a space with an action of an $E_\infty$ operad.
\end{definition}

In order to state the analog of Corollary \ref{23a} we need to use augmented
cosimplicial spaces. 

\begin{definition}
Let $X^\b$ be an augmented cosimplicial space.  Define
$\Tot(X^\b)$ to be 
\[
\Hom_{\Delta_+}(\Delta^\b,X^\b).
\]
\end{definition}

This can be described more simply:  $\Tot(X^\b)$ is the total space (in the
sense of Definition \ref{11aaaa}) of the restriction of $X^\b$ to $\Delta$.

\begin{theorem}
\label{32}
Let $X^\b$ be an augmented cosimplicial space with a map
\[
\langle f\rangle: X^{f^{-1}(1)}\times X^{f^{-1}(2)} \to X^T
\]
for each $f:T\to\{1,2\}$. 
Suppose that the maps $\langle f\rangle$ satisfy the analogs of
Propositions \ref{28}, \ref{29}, and \ref{30}, and that there is an element
$\varepsilon\in X^\emptyset$ satisfying the analog of Proposition \ref{31}.
Then $\Tot(X^\b)$ is an $E_\infty$ space.
\end{theorem}

\begin{remark}
We expect that $\Tot$ induces a Quillen equivalence between augmented
cosimplicial spaces satisfying the hypothesis of Theorem \ref{32} and
$E_\infty$ spaces.
\end{remark}

The remainder of this section gives an outline of the proof of Theorem
\ref{32}; for details see \cite{MS3}.

The proof follows the same pattern as the proof of Theorem \ref{17}.  The first
step is

\begin{proposition}
The category of augmented cosimplicial spaces has a symmetric monoidal 
structure $\boxtimes$ with the property that $X^\b$ satisfies the hypothesis 
of Theorem \ref{32} if and only if it is a commutative $\boxtimes$-monoid.
\end{proposition}

The basic idea in defining 
$X^\b \boxtimes Y^\b$ is that we build it from formal symbols 
$\langle f \rangle (x,y)$.  In order to get a cosimplicial object we have to
build in the cosimplicial operators, so we consider symbols of the form
\[
\phi_*(\langle f \rangle (x,y))
\]
where $f:T\to \{1,2\}$ and $\phi:T\to S$ is an order-preserving map: such a 
symbol will represent a point in the $S$-th space $(X^\b \boxtimes Y^\b)^S$.
We require these symbols to satisfy the relation in Proposition \ref{28}.

Our next two definitions make this precise.

\begin{definition}
Given $S\in
\Delta_+$, let $\I_S$ be
the category whose objects are diagrams 
\begin{equation}
\label{33}
\xymatrix{
\{1,2\}
&
T
\ar[l]_-f
\ar[r]^\phi
&
S
}
\end{equation}
where $T$ is a finite totally ordered set and $\phi$ is order-preserving, and
whose morphisms are commutative diagrams
\begin{equation}
\label{34}
\xymatrix{
\{1,2\}
\ar[d]_{=}
&
T
\ar[l]_-f
\ar[r]^\phi
\ar[d]_{\psi}
&
S
\ar[d]_{=}
\\
\{1,2\}
&
T'
\ar[l]_-{f'}
\ar[r]^{\phi'}
&
S
}
\end{equation}
with $\psi$ order-preserving.
\end{definition}

We will denote an object \eqref{33} of $\I_S$ by $(f,\phi)$. Given augmented
cosimplicial spaces $X^\b$ and $Y^\b$ we consider the functor from $\I_S$ to
spaces which takes $(f,\phi)$ to 
\[
X^{f^{-1}(1)}\times Y^{f^{-1}(2)}
\]
and a morphism \eqref{34} to the map
\[
(\psi_1)_* \times (\psi_2)_*
\]
where $\psi_i:f^{-1}(i)\to (f')^{-1}(i)$ is the restriction of $\psi$.

\begin{definition}
\label{34a}
Define $X^\b\boxtimes Y^\b$ by
\[
(X^\b\boxtimes Y^\b)^S=\colim_{(f,\phi)\in \I_S} X^{f^{-1}(1)}\times 
Y^{f^{-1}(2)}
\]
for $S\in \Delta_+$.
\end{definition}

The verification that $\boxtimes$ is a symmetric monoidal product is given in
\cite[Section 6]{MS3}.

\begin{remark}
Readers familiar with Kan extensions will recognize that $X^\b\boxtimes Y^\b$
is one; see \cite[Section 6]{MS3}.
\end{remark}

Next we apply Remark \ref{25}(b) with $D=\Delta^\b$ to get an operad
$\D$ whose $k$-th space is
\[
\Hom_\Delta(\Delta^\b,(\Delta^\b)^{\boxtimes k})
\]
If $X^\b$ is a commutative $\boxtimes$-monoid we define an action of $\D$ on
$\Tot(X^\b)$ by letting
\[
\theta:\D(k)\times (\Tot(X^\b))^k \to \Tot(X^\b)
\]
be the map that takes $(h,\tau_1,\ldots,\tau_k)$ to the composite
\[
\Delta^\b\xrightarrow{h}(\Delta^\b)^{\boxtimes k}
\xrightarrow{\tau_1\boxtimes\cdots\boxtimes \tau_k} (X^\b)^{\boxtimes k}
\xrightarrow{\mu} X^\b
\]
where $\mu$ is the monoidal structure map of $X^\b$.

To complete the proof of Theorem \ref{32} it only remains to show that each
$\D(k)$ is contractible; see \cite[Section 10]{MS3} for the proof of this.

\begin{remark}
\label{34b}
One can give a construction analogous to $\boxtimes$ for the category of
ordinary (unaugmented) cosimplicial spaces by requiring $f$ to be a surjection
in Definition \ref{34a}.  This gives a product which is coherently associative
and commutative but not unital.
\end{remark}

\section{The little $n$-cubes operad.}

\label{s9}

We have seen in Section \ref{s2} that $\Omega Z$ is an $A_\infty$ space.  For
$n\geq 2$ the space $\Omega^n Z$ has a commutativity property intermediate
between $A_\infty$ and $E_\infty$; moreover, $\Omega^n Z$ has stronger
commutativity than $\Omega^m Z$ if $n>m$. In this section we see how to make
this precise.

Fix $n\geq 1$.  Let $I$ denote the interval $[0,1]$. 

\begin{definition}
A TD-map $I^n\to I^n$ is a composite $T\circ D$, where $T$ is a translation and
$D$ is a dilation (i.e., multiplication by a scalar).
\end{definition}

A TD-map takes $(t_1,\ldots,t_n)$ to $(a_1+bt_1,\ldots,a_n+bt_n)$,
where $a_1,\ldots,a_n$ and $b$ are constants with $a_i\geq 0$, $b>0$ and
$a_i+b<1$.  The image of 
a TD-map is called a ``little $n$-cube.'' A TD-map is completely 
determined by its image.

\begin{definition}
(a) For $k\geq 0$, let $\C_n(k)$ be the space in which a point is a 
$k$-tuple $(\kappa_1,\cdots,\kappa_k)$ of TD-maps $I^n\to I^n$ such that the 
images of the $\kappa_i$ have disjoint interiors.  

(b) Let $\C_n$ be the collection of spaces $\{\C_n(k)\}_{k\geq 0}$.
\end{definition}

In the special case $n=2$, the elements of $\C_2(k)$ can be represented by 
pictures in the plane.  For example, the picture

\medskip

\noindent
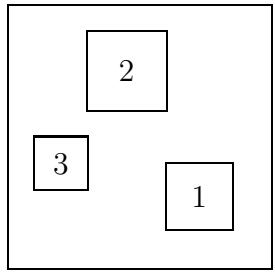
\mbox{}
\hfill
\begin{picture}(100,100)

\put(0,0){\line(1,0){100}}
\put(0,100){\line(1,0){100}}
\put(0,0){\line(0,1){100}}
\put(100,0){\line(0,1){100}}
\put(10,30){\framebox(20,20){3}} 
\put(30,60){\framebox(30,30){2}} 
\put(60,15){\framebox(25,25){1}} 

\end{picture}
\hfill
\mbox{}

\noindent
represents an element of $\C_2(3)$.

Next we define an operad structure for $\C_n$.  

\begin{definition}
\label{35}
(a)
Let $1\in \C_n(1)$ be the
identity map of $I^n$.  

(b)
Give $\C_n(k)$ the right $\Sigma_k$ action that
permutes the $\kappa_i$.  

(c)
Define 
\[
\gamma:\C_n(k)\times\C_n(j_1)\times\cdots\times \C_n(j_k)
\to \C_n(j_1+\cdots+j_k)
\]
as follows: if
\[
c=(\kappa_1,\cdots,\kappa_k)
\] 
is a point of $\C_n(k)$, and 
\[
d_i=(\lambda_{i1},\ldots,\lambda_{ij_i})
\]
is a point of $\C_n(j_i)$ for $1\leq i\leq k$, then 
$\gamma(c,d_1,\ldots,d_k)$ is the point 
$
(\nu_{11},\ldots,\nu_{kj_k}),
$
where
\[
\nu_{il}=\kappa_i\circ \lambda_{il}\,.
\]
\end{definition}

For example, if $n=2$ and $c\in \C_2(2)$, $d_1\in \C_2(3)$ and $d_2\in \C_2(2)$
are represented by

\medskip

\noindent
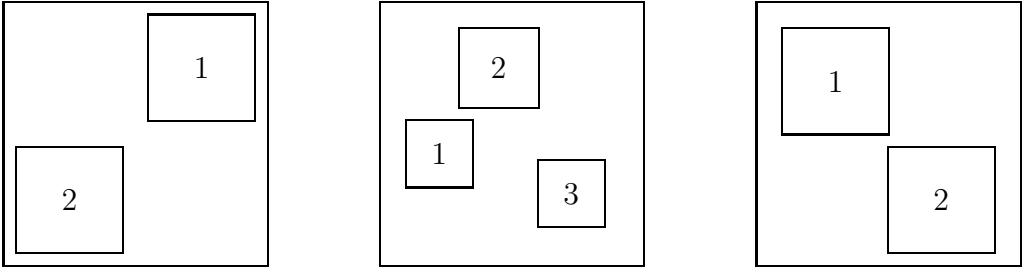
\mbox{}
\hfill
\begin{picture}(100,100)

\put(0,0){\line(1,0){100}}
\put(0,100){\line(1,0){100}}
\put(0,0){\line(0,1){100}}
\put(100,0){\line(0,1){100}}
\put(5,5){\framebox(40,40){2}} 
\put(55,55){\framebox(40,40){1}} 

\end{picture}
\hfill
\begin{picture}(100,100)

\put(0,0){\line(1,0){100}}
\put(0,100){\line(1,0){100}}
\put(0,0){\line(0,1){100}}
\put(100,0){\line(0,1){100}}
\put(10,30){\framebox(25,25){1}}
\put(30,60){\framebox(30,30){2}}
\put(60,15){\framebox(25,25){3}}

\end{picture}
\hfill
\begin{picture}(100,100)

\put(0,0){\line(1,0){100}}
\put(0,100){\line(1,0){100}}
\put(0,0){\line(0,1){100}}
\put(100,0){\line(0,1){100}}
\put(10,50){\framebox(40,40){1}}
\put(50,5){\framebox(40,40){2}}

\end{picture}
\hfill
\mbox{}

\medskip

\noindent 
respectively, then $\gamma(c,d_1,d_2)\in \C_2(5)$ is represented by

\medskip

\noindent
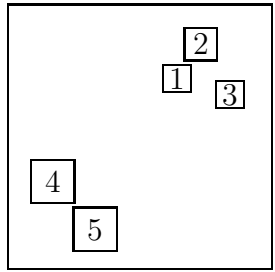
\mbox{}
\hfill
\begin{picture}(100,100)

\put(0,0){\line(1,0){100}}
\put(0,100){\line(1,0){100}}
\put(0,0){\line(0,1){100}}
\put(100,0){\line(0,1){100}}
\put(9,25){\framebox(16,16){4}}
\put(25,7){\framebox(16,16){5}}
\put(59,67){\framebox(10,10){1}}
\put(67,79){\framebox(12,12){2}}
\put(79,61){\framebox(10,10){3}}

\end{picture}
\hfill
\mbox{}

\medskip

\begin{remark}
\label{36}
(a) The definition of $\C_n$ and its composition maps is due to Boardman and 
Vogt \cite{BV1,BV2}. They were working in a somewhat different context (PROP's 
instead of operads).

(b) $\C_1$ is the operad generated by the non-symmetric operad $\A$ defined in
Section \ref{s2} (see Remark \ref{25}(c)).
\end{remark}

The reason for defining the operad $\C_n$ is that it acts on $\Omega^n Z$.
To describe this action we think of an element of $\Omega^n Z$ as a map
$I^n\to Z$ which takes the boundary of $I^n$ to the basepoint $*$ of $Z$.
Then
\[
\theta:\C_n(k)\times (\Omega^n Z)^k \to \Omega^n Z
\]
is defined as follows:
if 
\[
c=(\kappa_1,\ldots,\kappa_k)
\]
is a point of $\C_n(k)$ and $\alpha_i\in \Omega^n Z$
for $1\leq i\leq k$ then
$\theta(c,\alpha_1,\ldots,\alpha_k)$ is the map $I^n\to Z$
which is $\alpha_i\circ (\kappa_i)^{-1}$ on the image of $\kappa_i$ and $*$ for
points which are not in the image of any $\kappa_i$.  

For example, if $c$ is the element of $\C_2(3)$ represented by

\medskip

\noindent
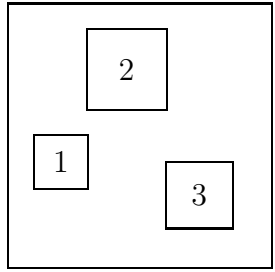
\mbox{}
\hfill
\begin{picture}(100,100)

\put(0,0){\line(1,0){100}}
\put(0,100){\line(1,0){100}}
\put(0,0){\line(0,1){100}}
\put(100,0){\line(0,1){100}}
\put(10,30){\framebox(20,20){1}} 
\put(30,60){\framebox(30,30){2}} 
\put(60,15){\framebox(25,25){3}} 

\end{picture}
\hfill
\mbox{}

\noindent
then $\theta(c,\alpha_1,\alpha_2,\alpha_3)$ is the map $I^2\to Z$ represented 
by the picture

\medskip

\noindent
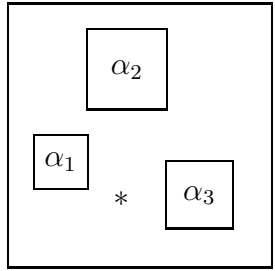
\mbox{}
\hfill
\begin{picture}(100,100)

\put(0,0){\line(1,0){100}}
\put(0,100){\line(1,0){100}}
\put(0,0){\line(0,1){100}}
\put(100,0){\line(0,1){100}}
\put(10,30){\framebox(20,20){$\alpha_1$}} 
\put(30,60){\framebox(30,30){$\alpha_2$}} 
\put(60,15){\framebox(25,25){$\alpha_3$}} 
\put(40,20){*}

\end{picture}
\hfill
\mbox{}

\medskip
\noindent
(where the $\alpha$'s in the picture are appropriately scaled).

As one would expect, there is an analog of Theorems \ref{10} and \ref{26}: if
$Y$ has a grouplike $\C_n$ action then $Y$ is weakly equivalent to $\Omega^n Z$
for some $Z$.  In fact something a bit more general is true; we pause to give
the relevant definitions, which will also be used in Section \ref{s10}.

\begin{definition} 
Let $\O$, $\O'$ be operads.  An {\it operad morphism} $\zeta:\O\to \O'$ is a
sequence of maps 
\[
\zeta_k:\O(k)\to \O'(k)
\]
such that 

(a) $\zeta_1$ takes the unit element in $\O(1)$ to that in $\O'(1)$,

(b) each $\zeta_k$ is $\Sigma_k$ equivariant, and

(c) the diagram
\[
\xymatrix{
\O(k)\times \O(j_1)\times \cdots \times \O(j_k)
\ar[d]_{\gamma}
\ar[rr]^-{\zeta_k\times\zeta_{j_1}\times \cdots}
&&
\O'(k)\times \O'(j_1)\times \cdots \times \O'(j_k)
\ar[d]_{\gamma}
\\
\O(j_1+\cdots+j_k)
\ar[rr]^-{\zeta_{j_1+\cdots+j_k}}
&&
\O'(j_1+\cdots+j_k)
}
\]
commutes for all $k,j_1,\ldots,j_k\geq 0$.
\end{definition}

\begin{definition}
A morphism $\zeta:\O\to\O'$ is a weak equivalence if each $\zeta_k$ is a weak
equivalence of spaces.  Two operads $\O$ and $\O'$ are weakly equivalent if
there is a diagram of operads and weak equivalences of operads
\[
\O \leftarrow \cdots \rightarrow \O'
\]
\end{definition}

\begin{definition}
An operad is an $E_n$ operad if it is weakly equivalent to $\C_n$.
\end{definition}

Now the analog of Theorems \ref{10} and \ref{26} is

\begin{theorem}
\label{37}
$Y$ is weakly
equivalent to $\Omega^n Z$ for some space $Z$ $\Longleftrightarrow$
$Y$ has a grouplike action of an $E_n$ operad. 
\end{theorem}

The $\Longrightarrow$ direction, and the $\Longleftarrow$ direction for
connected $Y$, are due to Boardman and Vogt 
\cite{BV1,BV2}. 
A simpler proof of the $\Longleftarrow$ direction for connected $Y$ was given 
by May in \cite{MayG}.  The general 
case of the $\Longleftarrow$ direction is due to May \cite{MayE}.

Theorem \ref{37} is aesthetically pleasing but has not often been applied
because it is usually hard to show that a space is an $E_n$ space (for 
$1<n<\infty$) without knowing in advance that it is an $n$-fold loop space.  
In Section \ref{s10} we will address this difficulty by giving a sufficient 
condition for Tot of a cosimplicial space to be an $E_n$ space.

Since this section is intended as an introduction to $\C_n$, we should mention
that the most important uses of $\C_n$ in algebraic topology come from the
``approximation theorem'' \cite[Theorem 2.7]{MayG}.  This theorem gives a model
for $\Omega^n \Sigma ^n Z$, built from $Z$ and $\C_n$.  Using this model, Fred
Cohen has given a complete description of the homology of $\Omega^n \Sigma ^n
Z$ \cite{CLM}.  Another basic fact is that the model splits stably as a
wedge of pieces of the form
\[
\C_n(k)_+\wedge_{\Sigma_k} Z^{\wedge k}
\]
(where $+$ means add a disjoint basepoint and $\wedge$ is the smash product);
this is called the Snaith splitting \cite{Snaith,RCohen}.

\section{A sufficient condition for $\Tot(X^\b)$ to be an $E_n$ space.}

\label{s10}

In this section we give the analog of Theorem \ref{32} for $E_n$ actions.  
The hypothesis of Theorem \ref{32} refers to $\langle f\rangle$ operations
where $f$ ranges through functions $T\to\{1,2\}$. 
For our current purpose we need to consider functions $f:T\to\{1,\ldots,k\}$
for {\it all} $k\geq 2$ but we will only use those of
``complexity $\leq n$'' (see Definitions \ref{38a} and \ref{38b}).

Let us return to the motivating example $S^\b W$.  The extension of the
definition of $\langle f\rangle$ to functions $f:T\to \{1,\ldots,k\}$ is
routine:

\begin{definition}
Given $f:T\to\{1,\ldots,k\}$ define a natural transformation
\[
\langle f \rangle: S^{f^{-1}(1)} W\times \cdots\times S^{f^{-1}(k)} W
\to S^T W
\]
by the equation
\[
\langle f \rangle(x_1,\ldots,x_k)(\sigma)=
x_1(\sigma(f^{-1}(1)))\cdot \,\cdots \,\cdot x_k(\sigma(f^{-1}(k)))
\]
for $\sigma\in S_{T} W$; here $\cdot$ is multiplication in $\mathbb Z$.
\end{definition}

These operations satisfy properties analogous to Propositions \ref{28},
\ref{29}, \ref{30} and \ref{31}; the precise formulation is left to the reader
(see \cite[Definitions 9.3--9.7]{MS3} for a hint).

\begin{remark}
\label{38}
$\langle f\rangle$ operations with $k>2$ are composites of those with $k=2$; 
that is why we were
able to restrict to $k=2$ in Sections \ref{s7} and \ref{s8}.  However, it is
not true that an operation $\langle f\rangle$ with complexity $\leq n$ (see
Definitions \ref{38a} and \ref{38b}) can be decomposed into operations with 
$k=2$ {\it and} complexity $\leq n$, which is why we cannot restrict to $k=2$ 
in this section.
\end{remark}

As we have seen in Section \ref{s9}, an $E_n$ operad encodes commutativity
which is intermediate between $A_\infty$ (no commutativity) and $E_\infty$
(full commutativity).  We therefore want, for each $n$, a family of operations
which interpolates between $\sqcup$ and the family of all $\langle f \rangle$ 
operations.  Notice that, in general, the ordered sets 
$f^{-1}(1),\ldots,f^{-1}(k)$ are mixed together in $T$, but when $f$
corresponds to an iterate of $\sqcup$ they are not mixed. We therefore
introduce a way of measuring the amount of mixing.

We begin with the case $k=2$.
First observe that a function $f$ from a finite totally ordered set to
$\{1,2\}$ is the same thing as a finite sequence of 1's and 2's.

\begin{definition}
\label{38a}
(a)
The complexity of a finite sequence of 1's and 2's is the number
of times the sequence changes from 1 to 2 or from 2 to 1.

(b)
The complexity of a function $f:T\to \{1,2\}$ is the complexity of the 
corresponding sequence.
\end{definition}

For example, if $f$ corresponds to the sequence 11222122112 then the
complexity of $f$ is 5.

Next let $k>2$.  A function $f:T\to \{1,\ldots,k\}$ corresponds to a finite
sequence with values in $\{1,\ldots,k\}$.  We consider the subsequences that
have only two values: for example in the sequence 12313212 we consider the
subsequences 121212, 23322 and 13131.  As in Definition \ref{38a}, the
complexity of such a subsequence is the number of times it changes its value.

\begin{definition}  
\label{38b}
(a) The complexity of a sequence with values in $\{1,\ldots,k\}$ is the maximum
of the complexities of the subsequences with only two values.

(b) The complexity of $f:T\to\{1,\ldots,k\}$ is the complexity of the
sequence corresponding to $f$.
\end{definition}  

In the example just given, the complexity of 121212 is 5, the complexity of
23322 is 2, and the complexity of 13131 is 4, so the complexity of 12313212
is 5.

\begin{remark}
The definition of complexity is suggested by \cite{Smith};
the reason we use it here is that it is well-adapted to the proof of Theorem
\ref{39} below.
There may be other ways of defining complexity that would also lead to
Theorem \ref{39}, although this seems unlikely.
\end{remark}

We can now state the analog of Theorem \ref{32}.

\begin{theorem}
\label{39}
Fix $n$.
Let $X^\b$ be an augmented cosimplicial space with a map
\[
\langle f\rangle: X^{f^{-1}(1)}\times \cdots \times X^{f^{-1}(k)} \to X^T
\]
for each $f:T\to\{1,\ldots,k\}$ with complexity $\leq n$.
Suppose that the maps $\langle f\rangle$ are consistent with the cosimplicial
operators (in the sense of \cite[Definition 9.4]{MS3}) and are commutative, 
associative and unital (in the sense of \cite[Definitions 9.5, 9.6 and 
9.7]{MS3}).  Then $\Tot(X^\b)$ is an $E_n$ space.
\end{theorem}

We expect that $\Tot$ induces a Quillen equivalence between augmented
cosimplicial spaces satisfying the hypothesis of Theorem \ref{39} and
$E_n$ spaces.

The proof of Theorem \ref{39} is similar in outline to the proofs of
Theorems \ref{17} and \ref{32}.   However, just as $E_n$ interpolates between
$A_\infty$ and $E_\infty$, we need a way to interpolate between the concepts of
monoidal product and symmetric monoidal product.  The next section is
devoted to this.

\section{An extension of Remark \ref{25}(b).}

\label{s11}

Remark \ref{25}(b) says that a symmetric monoidal 
product $\boxtimes$ on a topological category $\U$, together with a choice of 
an object $D\in \U$, leads to an operad $\O$.  We begin with an outline of the
proof of this fact.

$\boxtimes$ is a binary operation, so the first step
is to choose, for each $k>2$, a specific way of inserting
parentheses\footnote{A different choice gives a naturally isomorphic 
functor, by MacLane's coherence theorem \cite[Section VII.7]{MacLane}.}
to get a $k$-fold iterate of $\boxtimes$ which we denote by
\[
\boxtimes^k: \U^k\to \U.
\]
We also define $\boxtimes^1$ to be the identity functor and $\boxtimes^0$ to be
the unit object of $\boxtimes$.

Next we define the spaces of the operad $\O$ by
\[
\O(k)=\Hom_\U(D,\boxtimes^k(D,\ldots,D))
\]
for $k\geq 0$.

To define the action of $\Sigma_k$ on $\O(k)$, we use MacLane's 
coherence theorem.  This gives, for each $\sigma\in \Sigma_k$, a natural 
isomorphism 
\[
\sigma_*:\boxtimes^k(X_1,\ldots,X_k)
\to
\boxtimes^k(X_{\sigma(1)},\ldots,X_{\sigma(k)})
\]
and in particular a self-map of $\boxtimes^k(D,\ldots,D)$.  

We use the coherence theorem again
to get a natural isomorphism
\[
\Gamma:\boxtimes^k(\boxtimes^{j_1},\ldots,\boxtimes^{j_k})
\to
\boxtimes^{j_1+\cdots+j_k}
\]
which induces the structure map $\gamma$ of $\O$.  

It remains to check that $\gamma$ has the associativity, unitality and 
equivariance properties required by the definition of operad; for this we 
apply the coherence theorem one more time to see that $\Gamma$ has 
associativity, unitality and equivariance properties (see 
\cite[Section 4]{MS3} for the explicit statements) from which those for 
$\gamma$ can be deduced.  This completes the proof of Remark \ref{25}(b).

\bigskip

The same proof proves something more general. Assume that for each $k$ we are 
given a subfunctor of $\boxtimes^k$, that is, a functor $\F_k$ with a natural 
monomorphism to $\boxtimes^k$.  Assume further that $\F_k$ is closed under 
$\sigma_*$ (that is, $\sigma_*$ takes $\F_k(X_1,\ldots,X_k)$ 
to $\F_k(X_{\sigma(1)},\ldots,X_{\sigma(k)})$) and that the 
collection $\{\F_k\}_{k\geq 0}$ is closed under $\Gamma$.  The argument given 
above shows:

\begin{proposition}
\label{40}
Under these assumptions the collection $\{\Hom_\U(D,\F_k(D,\ldots,D))\}_{k\geq
0}$ is an operad, with $\Sigma_k$ action induced by the maps $\sigma_*$ and 
$\gamma$ induced by $\Gamma$.
\end{proposition}

\begin{remark}
In \cite[Section 4]{MS3} we give a more general version of \ref{40}, using 
the concept of
``functor-operad.''
A functor-operad is a collection of functors
\[
\F_k:\U^k\to \U
\]
with just enough structure to satisfy the conclusion of Proposition \ref{40}.
This concept has been discovered independently, in a different context and in 
a more general form, by Batanin \cite{Bat02}, who calls them ``internal 
operads.''
\end{remark}

\section{Proof of Theorem \ref{39}.} 

\label{s12}

Recall that in the proof of Theorem \ref{32} we constructed a symmetric 
monoidal product $\boxtimes$ on the category of augmented cosimplicial 
spaces.  Our first task is to give a formula for the iterate $\boxtimes^k$. 

Given $S\in \Delta_+$ and $k\geq 0$ let $\I_S(k)$ be the category whose objects
are diagrams
\begin{equation}
\label{41}
\xymatrix{
\{1,\ldots,k\}
&
T
\ar[l]_-f
\ar[r]^\phi
&
S
}
\end{equation}
where $T$ is a finite totally ordered set and $\phi$ is order-preserving, and
whose morphisms are commutative diagrams
\[
\xymatrix{
\{1,\ldots,k\}
\ar[d]_{=}
&
T
\ar[l]_-f
\ar[r]^\phi
\ar[d]_{\psi}
&
S
\ar[d]_{=}
\\
\{1,\ldots,k\}
&
T'
\ar[l]_-{f'}
\ar[r]^{\phi'}
&
S
}
\]
with $\psi$ order-preserving.

We will denote an object \eqref{41} of $\I_S$ by $(f,\phi)$.

\begin{definition}
Let $X_1^\b,\ldots,X_k^\b$ be augmented cosimplicial spaces.
Define $\Xi_k(X_1^\b,\ldots,X_k^\b)$ to be the cosimplicial space whose value
at $S\in \Delta_+$ is 
\[
\colim_{(f,\phi)\in\I_S(k)} X_1^{f^{-1}(1)}\times \cdots\times X_k^{f^{-1}(k)}
\]
\end{definition}

When $k=2$ this is the same as the definition of $X_1^\b\boxtimes X_2^\b$.  
In general we have

\begin{proposition}
$\boxtimes^k$ is naturally isomorphic to $\Xi_k$. 
\end{proposition}

For the proof see \cite[Section 6]{MS3} (also cf.\ Remark \ref{38}).

Now fix $n$, and let $\I_S^n(k)$ be the full subcategory of $\I_S(k)$ whose
objects are the $(f,\phi)$ for which the complexity of $f$ is $\leq n$.

\begin{definition}
Let $X_1^\b,\ldots,X_k^\b$ be augmented cosimplicial spaces.
Define $\Xi_k^n(X_1^\b,\ldots,X_k^\b)$ to be the cosimplicial space whose value
at $S\in \Delta_+$ is
\[
\colim_{(f,\phi)\in\I_S^n(k)} X^{f^{-1}(1)}\times \cdots\times X^{f^{-1}(k)}
\]
\end{definition}

\begin{proposition}
For each $n$,
the sequence of functors $\{\Xi^n_k\}_{k\geq 0}$ satisfies the hypothesis of 
Proposition \ref{40}.
\end{proposition}

For the proof see \cite[Section 8]{MS3}.

Applying Proposition \ref{40} we obtain an operad with $k$-th space
\begin{equation}
\label{42}
\Hom_{\Delta_+}(\Delta^\b,\Xi^n(\Delta^\b,\ldots,\Delta^\b))
\end{equation}
We will denote this operad by $\D_n$.

If $X^\b$ satisfies the hypothesis of Theorem \ref{39} then there are maps
\[
\xi_k:\Xi_k^n(X^\b,\ldots,X^\b)\to X^\b
\]
for each $k\geq 0$. We define an action of $\D_n$ on $\Tot(X^\b)$ by letting
\[
\theta:\D_n(k)\times (\Tot(X^\b))^k \to \Tot(X^\b)
\]
be the map that takes $(h,\tau_1,\ldots,\tau_k)$ to the composite
\[
\Delta^\b\xrightarrow{h}\Xi_k^n(\Delta^\b,\ldots,\Delta^\b)
\xrightarrow{\Xi_k^n(\tau_1,\cdots,\tau_k)} 
\Xi_k^n(X^\b,\ldots,X^\b)
\xrightarrow{\xi_k} X^\b
\]

To complete the proof of Theorem \ref{39} it remains to show that $\D_n$
is weakly equivalent to $\C_n$.  This is more difficult than the corresponding
step in the proofs of Theorems \ref{17} and \ref{32} because in those cases it
was only necessary to show that certain spaces were contractible, whereas 
here we need to show not just that the spaces $\D_n(k)$ and $\C_n(k)$ are
weakly equivalent but that the operad structures are compatible.
The proof is given in \cite[Section 12]{MS3}; the basic idea is to show that
the operads $\C_n$ and $\D_n$ can be written as homotopy colimits, over the
same indexing category, of contractible sub-operads. 

\begin{remark}
\label{rev2}
The operad $\D_n$ is of interest in its own right, as an $E_n$ operad whose
structure is in some ways simpler than that of $\C_n$.  In \cite[Section
11]{MS3} it is shown that $\D_n(k)$ is homeomorphic to 
\[
Z_k^n\times \Tot(\Delta^\b)
\]
where $Z_k^n$ is the zeroth space of $\Xi_k^n(\Delta^\b,\ldots,\Delta^\b)$.
Moreover, the space $Z_k^n$ has an explicit cell decomposition which is
well-related to the operad structure of $\D_n$.  The cellular chain complexes
of the $Z_k^n$ form a chain operad which is studied in \cite{MS2} (where it is
called $\S_n$).
\end{remark}

\section{Applications.}

\label{s13}

\subsection{The topology of a space of knots.}

The space of imbeddings of $S^1$ in ${\mathbb R}^k$ is of considerable interest
in topology.  It turns out that a closely related space satisfies the
hypothesis of Theorem \ref{39} with $n=2$, and is therefore a two-fold loop
space.

To be specific, let us consider the manifold-with-boundary 
${\mathbb R}^{k-1}\times I$. Fix points $x_0\in {\mathbb R}^{k-1}\times 
\{0\}$ and $x_1\in {\mathbb R}^{k-1}\times \{1\}$, and also fix tangent 
vectors $v_0$ at $x_0$ and $v_1$ at $x_1$.
Let Emb$(I,{\mathbb R}^{k-1}\times I)$ be
the space of embeddings of $I$ in ${\mathbb R}^{k-1}\times I$ which take
$0$ and $1$ to $x_0$ and $x_1$ respectively, with tangent vectors $v_0$ at 0
and $v_1$ at 1.  Let Imm$(I,{\mathbb R}^{k-1}\times I)$ be the 
analogous space with immersions instead of imbeddings.  Finally, let 
Fib$(I,{\mathbb R}^{k-1}\times I)$ be the fiber of the forgetful map
\[
\text{Emb}(I,{\mathbb R}^{k-1}\times I) \to 
\text{Imm}(I,{\mathbb R}^{k-1}\times I) 
\]
It follows from a theorem of Hirsch and Smale that 
Imm$(I,{\mathbb R}^{k-1}\times I)$ is homotopy equivalent to $\Omega 
S^{k-1}$, so Fib$(I,{\mathbb R}^{k-1}\times I)$ contains most of the 
information in Emb$(I,{\mathbb R}^{k-1}\times I)$.

Now assume $k\geq 4$.

Dev Sinha
\cite{S} (building on earlier work of
Goodwillie and Weiss) has given 
a cosimplicial space $X^\b$ with $\Tot(X^\b)$ weakly equivalent 
to Fib$(I,{\mathbb R}^{k-1}\times I)$. 
He has also shown that $X^\b$ satisfies the 
hypothesis of Theorem \ref{39} with $n=2$.  It follows that 
Fib$(I,{\mathbb R}^{k-1}\times I)$ is a two-fold loop space. 

\begin{remark}
When a cosimplicial space $X^\b$ satisfies the hypothesis of Theorem \ref{39}, 
the spectral sequence converging to the homology of $\Tot(X^\b)$ will have
extra structure coming from the $\langle f\rangle$ operations.  This should be
useful for analyzing the spectral sequence that converges to the homology of
Fib$(I,{\mathbb R}^{k-1}\times I)$.
\end{remark}

\subsection{Topological Hochschild Cohomology.}

Theorems \ref{17}, \ref{32} and \ref{39} are still true, with essentially the
same proofs, for cosimplicial
spectra (except that Cartesian products in the category of spaces are
replaced by smash products in the category of spectra).

The definition of Hochschild cohomology for associative rings (which will be
recalled in Section \ref{s16}) has an analog for associative ring spectra in 
the sense of \cite{EKMM}  or \cite{HSS}. 
If $R$ is an associative ring 
spectrum there is a cosimplicial spectrum $TH^\b(R)$ 
(see \cite[Example 3.4]{MS1} for the definition)
whose total spectrum 
$\Tot(TH^\b(R))$ is called the topological Hochschild cohomology spectrum of 
$R$.

In \cite{MS1} it is
shown that $TH^\b(R)$ satisfies the hypothesis of Theorem \ref{39} with $n=2$,
and therefore the topological Hochschild cohomology spectrum of $R$ is an $E_2$
spectrum.  This is a spectrum analog of Deligne's Hochschild cohomology
conjecture (see Section \ref{s16}).

\section{The framed little disks operad.}

\label{s14}

The framed little disks operad 
was defined by Getzler in \cite{Getzler};
it is a variant of the little 2-cubes operad. 

Let $B$ denote the closed unit disk in ${\mathbb R}^2$.

\begin{definition}
A TDR-map $B\to B$ is a composite $T\circ D\circ R$, where $T$ is a 
translation, $D$ is a dilation and $R$ is a rotation.
\end{definition}

\begin{definition}
\label{43}
(a) For $k\geq 0$, let $\F(k)$ be the space in which a point is a
$k$-tuple $(\kappa_1,\cdots,\kappa_k)$ of TDR-maps $B\to B$ such that the
images of the $\kappa_i$ have disjoint interiors.

(b) Let $\F$ be the collection of spaces $\{\F(k)\}_{k\geq 0}$.
\end{definition}

$\F$ is an operad, where the $\Sigma_k$ action on $\F(k)$ permutes the
$\kappa_i$ and the definition of $\gamma$ is analogous to Definition 
\ref{35}(c).

\begin{remark}
It is instructive to consider the relation between $\F$ and $\C_2$.

(a)
If we require the $\kappa_i$ in Definition \ref{43} to
be $TD$ maps (that is, composites of translations and dilations), we get a
suboperad $\F_0$ of $\F$.  
By restricting TD maps $B\to B$ to the square 
inscribed in $B$ we get an equivalence of operads $\F_0 \to \C_2$.  

(b) The $k$-th space of $\F$ is the Cartesian product
$\F_0(k)\times (S^1)^k$ (but note that the projections $\F(k) \to \F_0(k)$ do 
not give a map of operads).  

(c)
An action of $\F$ on a space $X$ is the same thing as an $\F_0$ action 
together with a suitably compatible $S^1$ action.
\end{remark}

\begin{remark}
One reason that $\F$ is important is that an $\F$ action on a space $X$ induces
a Batalin-Vilkovisky structure on $H_*X$ (see \cite{Getzler}).
\end{remark}

The analog of Theorem \ref{39} for $\F$ actions has a surprisingly simple
form.  As motivation we consider the following situation:  let $V_\b$ be a
cyclic set, that is, a simplicial set together with maps 
\[
t:V_m\to V_m
\]
for each $m\geq 0$
satisfying certain relations with the simplicial operators (see
\cite[Definition 9.6.1]{Weibel}).  Define $A^\b$ to be
$\Map(V_\b,{\mathbb Z})$.  Then 
$A^\b$ is a cocyclic abelian group,
that is, it is a cosimplicial abelian group together with maps 
\[
\tau:A^m\to A^m
\]
for each $m\geq 0$ which satisfy appropriate relations with the cosimplicial
operators.  We can define $\sqcup$ and the other
$\langle f \rangle$ operations on $A^\b$ in analogy with Definition \ref{26a}
(the precise definition is left as an exercise for the reader; the basic idea
is to use iterated face maps to interpret the symbol $\sigma(U)$ in this
context). The relations between the maps $t$ and the simplicial operators
imply that all $\langle f\rangle$ operations of
complexity $\leq 2$ are generated by $\sqcup$ and $\tau$, subject to the 
relation
\begin{equation}
\label{44}
\tau^{p+1}(x\sqcup y)=y\sqcup x,
\end{equation}
where $x$ is in $A^p$ and $\tau^{p+1}$ denotes the $(p+1)$-st iterate of $\tau$.

\begin{theorem}
\label{S2}
If $X^\b$ is a cocyclic space with a product
\[
\sqcup:X^p\times X^q \to X^{p+q+1}
\]
which is associative and unital
and satisfies
\eqref{18}, \eqref{19} and \eqref{44}
then $\Tot(X^\b)$ has an action of $\F$.
\end{theorem}

The proof is similar to that of Theorem \ref{39}; see \cite{MS5}.

\section{Cosimplicial chain complexes.}

\label{s15}

The theory developed in Sections \ref{s4}, \ref{s5}, \ref{s8},
\ref{s10}, \ref{s12} and \ref{s14} 
has a precise analog with spaces replaced by chain complexes; see 
\cite{MS4}.  In this section we give a brief discussion.

First we need the appropriate concept of operad.  In fact one can define
nonsymmetric operads in any monoidal category by replacing the Cartesian 
products in Definition \ref{6} by the monoidal product, and one can define
operads in any symmetric monoidal category by analogy with Definition 
\ref{24a}.  The category of chain complexes is a symmetric
monoidal category (the monoidal product is the usual tensor product of
chain complexes) and operads in this category are called {\it
chain operads.}

\begin{definition}
(a)
A chain complex is {\it weakly contractible} if its homology is ${\mathbb Z}$
in dimension 0 and zero in all other dimensions.  

(b)
An $A_\infty$ chain operad is
a nonsymmetric chain operad $\O$ for which each $\O(k)$ is a weakly
contractible chain complex.  

(c)
An $E_\infty$ chain operad is a chain operad $\O$
for which each $\O(k)$ is a weakly contractible chain complex.\footnote{For
technical reasons, it is usual to require
in addition that the action of each $\Sigma_k$ should be free. The operad
$\T$ defined below has this property.}
\end{definition}

Next we need the analog of Tot.
We have already defined the conormalization of a cosimplicial abelian group
(Definition \ref{11a}).
We now extend this definition to cosimplicial chain complexes.
Recall the cosimplicial chain complex $\Delta^\b_*$ (Definition \ref{S1}).

\begin{definition}
Let $B^\b_*$ be a cosimplicial chain complex.  The {\it conormalization} of
$B_*^\b$, denoted ${\mathbf C}(B_*^\b)$, is the cochain complex
\[
\Hom_\Delta (\Delta^\b_*,B_*^\b) \subset
\prod_{m=0}^\infty \Hom(\Delta^m_*,B_*^m).
\]
Here $\Hom_\Delta$ is Hom in the category of cosimplicial graded abelian
groups and the differential
is induced by the differentials of $\Delta^\b_*$ and $B^\b_*$.
\end{definition}

In practice it's useful to have an elementary description of ${\mathbf 
C}(B_*^\b)$.
First note that by fixing the internal degree $m$ we get a cosimplicial 
abelian group $B_m^\b$
and hence a cochain complex ${\mathbf C}(B_m^\b)$ (see Remark \ref{12} for 
elementary descriptions of this cochain complex).  The differential in 
$B_*^\b$ induces a differential
\[
{\mathbf C}(B_m^\b) \to {\mathbf C}(B_{m-1}^\b)
\]
so the cochain complexes ${\mathbf C}(B_m^\b)$ assemble into a bicomplex
(with differentials lowering degree in one direction and raising degree in the
other).  The conormalization of $B_*^\b$ is the 
totalization of this bicomplex: 
\[
{\mathbf C}(B_*^\b)^p=\prod_m {\mathbf C}(B^\b_m)^{p+m};
\]
in general this is an infinite product.

Next we observe that the definition of $\boxtimes$ in Section \ref{s8} has an
analog for augmented cosimplicial chain complexes, with $\times$ replaced by
$\otimes$.  As a consequence we get a chain operad $\T$ with 
\[
\T(k)={\mathbf C}((\Delta_*^\b)^{\boxtimes k})
\]

\begin{theorem}
(a) $\T$ is an $E_\infty$ chain operad.

(b) If $B_*^\b$ is a cosimplicial chain complex satisfying the hypothesis of
Theorem \ref{32} then $\T$ acts on ${\mathbf C}(B_*^\b)$.
\end{theorem}

\begin{remark}
\label{rev4}
The definition of $\T$ looks complicated, but in fact $\T$ has a simple 
explicit description.  For each fixed $k\geq 1$ and $q,r\geq 0$ let 
$\U_{q,r}(k)$ be the free abelian group generated by the symbols
\[
\xymatrix{
\{1,\ldots,k\}
&
[q]
\ar[l]_-f
\ar[r]^\phi
&
[r]
}
\]
where 
\begin{description}
\item[\rm (a)] $f$ is onto,
\item[\rm (b)] the image of $\phi$ contains all of $\{1,\ldots,r\}$ (but is
allowed to not contain 0),
\item[\rm (c)] $\phi$ is order-preserving, 
\item[\rm (d)] $\phi(i)=\phi(i+1)\Rightarrow f(i)\neq f(i+1)$.
\end{description}
Then the $p$-th group of the chain complex $\T(k)$ is
\[
\prod_q \U_{q,q+1-p-k}(k)
\]
It is not hard to show this from the definitions; see \cite{MS4}. 
\cite{MS4} also gives explicit formulas for the differential of $\T(k)$ and for
the operad structure maps of $\T$.
\end{remark}

\begin{remark}
\label{rev3}
$\T$ is not the same as the $E_\infty$ chain operad $\S$ defined in 
\cite{MS2}, but they are related: $\S$ can be obtained from $\T$ by the
condensation process described in \cite[Section 7]{MS1}.  We will show in
\cite{MS4} that the structural formulas for $\S$ can be deduced from 
those for $\T$; this is less elementary than the treatment of the structural
formulas in \cite{MS2} but avoids the eight pages of sign verifications in 
that paper.  Each of $\S$ and $\T$ has advantages: $\S$ is much smaller but 
$\T$ has useful formal properties (see Section \ref{s16c}).
\end{remark}

Next we need the definition of weak equivalence for chain operads.

\begin{definition}
A morphism $\zeta:\O\to\O'$ of chain operads is a weak equivalence if each 
$\zeta_k$ is a homology isomorphism.
Two chain operads $\O$ and $\O'$ are weakly equivalent if
there is a diagram of operads and weak equivalences of operads
\[
\O \leftarrow \cdots \rightarrow \O'
\]
\end{definition}

Now fix $n\geq 1$.  Applying the normalized singular chain functor to the
little $n$-cubes operad $\C_n$
we obtain a chain operad $S_* \C_n$.

\begin{definition}
An $E_n$ chain operad is a chain operad weakly equivalent to $S_*\C_n$.
\end{definition}

The definition of $\Xi_k^n$ in Section \ref{s12} has an analog for 
augmented cosimplicial chain complexes.  As a consequence we get a chain operad
$\T_n$ with
\[
\T_n(k)={\mathbf C}(\Xi^n_k(\Delta_*^\b,\ldots,\Delta_*^\b))
\]

\begin{theorem}
(a) $\T_n$ is an $E_n$ chain operad.

(b) If $B_*^\b$ is a cosimplicial chain complex satisfying the hypothesis of
Theorem \ref{39} then $\T_n$ acts on ${\mathbf C}(B_*^\b)$.
\end{theorem}

\begin{remark}
$\T_n$ has an explicit description similar to that in Remark \ref{rev4}, except
that $f$ is required to have complexity $\leq n$; see \cite{MS4}.  Also in 
\cite{MS4}, we show that the chain operad $\S_n$ defined in \cite{MS2} can be 
obtained from $\T_n$ by condensation.
\end{remark}

\begin{remark}
\label{S3}
The theory described in Section \ref{s14} also has a chain analog.  In
\cite{MS5} we construct a chain operad $\G$ which is weakly equivalent to
$S_*\F$, and we show that if $B_*^\b$ is a cosimplicial chain complex
satisfying the hypothesis of
Theorem \ref{S2} then $\G$ acts on ${\mathbf C}(B_*^\b)$.
\end{remark}

\section{Applications.}

\label{s16}

\subsection{Deligne's Hochschild cohomology conjecture.}

\label{s16a}

Let $R$ be an associative ring.  The {\it Hochschild cochain complex}\/ 
$C^*(R)$ is the cochain complex which in degree $p$ is 
\[
\Hom_{\mathbb Z}(R^{\otimes p}, R);
\]
the differential is determined by the formula
\begin{multline*}
(d\rho)(r_1\otimes\cdots\otimes r_{p+1}) \\
=r_1\rho(r_2\otimes\cdots)
+\sum_{i=1}^p (-1)^i \rho(\cdots \otimes r_i r_{i+1} \otimes\cdots)
+(-1)^{p+1} \rho(\cdots\otimes r_p)r_{p+1}
\end{multline*}
where $\rho \in C^p(R)$.  The {\it Hochschild cohomology}\/ $H^*(R)$ is the
cohomology of this complex.

Hochschild defined a cup product on $C^*(R)$ by 
\[
(\rho_1 \smallsmile \rho_2)(r_1\otimes\cdots\otimes r_{p+q})
=\rho_1(\cdots\otimes r_p)\cdot \rho_2(r_{p+1}\otimes\cdots)
\]
where $\rho_1\in C^p(R)$ and $\rho_2\in C^q(R)$.  This induces a product, also
denoted by $\smallsmile$, on $H^*(R)$.

Gerstenhaber showed in 1963 (see \cite{Gerst}) that $H^*(R)$ is what is now 
known as a {\it Gerstenhaber algebra}.  That is, he showed that the cup product 
on $H^*(R)$ is commutative and that there is a Lie bracket 
\[
[\ ,\ ]:H^p(R)\otimes H^q(R)\to H^{p+q-1}(R)
\]
such that $[x,\ ]$ is a derivation with respect to $\smallsmile$ for each $x\in
H^*(R)$.

About 10 years later, Fred Cohen showed that if $X$ has a $\C_2$ action then
$H_*X$ is a Gerstenhaber algebra (but he didn't use this terminology); see
\cite{CLM}.

In 1993, Deligne asked in a letter \cite{Deligne} whether these two 
examples of Gerstenhaber algebras were related: specifically he asked whether 
the cup product and Lie bracket on $H^*(R)$ are induced by an action of an 
$E_2$ chain operad on $C^*(R)$.

One reason this conjecture is important is because of its connection with
Kontsevich's deformation quantization theorem; see \cite{K1}.

The conjecture has been proved by several authors using quite different
methods (see \cite{T1,T2,MS1,V,K1,K2,HKV,Kauf}).  In \cite[Section 2]{MS2} we 
gave an especially simple proof by showing that the $E_2$ operad $\S_2$ 
defined in that paper acts on $C^*(R)$ by explicit formulas.  In \cite{MS4} 
we show that the $E_2$ operad $\T_2$ (see Section \ref{s15}) acts on 
$C^*(R)$, also by explicit formulas; this argument has the advantage that it 
avoids the complicated sign verifications needed in \cite{MS2}.

\subsection{Strong Frobenius algebras.}

\label{s16b}

By a {\it strong Frobenius algebra}\/ we mean an algebra $A$ over a field such
that $A$ is isomorphic to $A^*$ as an $A$-bimodule.  In \cite{MS5} we show
that if $A$ is a strong Frobenius algebra then the chain operad $\G$ (see
Remark \ref{S3}) acts on $C^*(A)$.  This is a strong form of Deligne's
Hochschild cohomology conjecture for these algebras.

\subsection{A theorem of Kriz and May.}

\label{s16c}

Let $\text{Ab}^\Delta$ denote the category of cosimplicial abelian groups, and
let $\text{Ch}^*_{\geq 0}$ denote the category of non-negatively graded cochain
complexes.

The conormalization functor gives an equivalence of categories
\begin{equation}
\label{S4}
{\mathbf C}:\text{Ab}^\Delta \to \text{Ch}^*_{\geq 0}
\end{equation}
(cf.\ \cite[Section 8.4]{Weibel}).

We mentioned in Remark \ref{34b} that the category of cosimplicial spaces has a
non-unital symmetric monoidal product $\boxtimes$; essentially the same
construction (with $\times$ replaced by $\otimes$) gives a non-unital symmetric
monoidal product $\boxtimes$ for $\text{Ab}^\Delta$.  It is natural to ask how
$\boxtimes$ is related to the equivalence \eqref{S4}, and this question has a
simple answer:

\begin{theorem}
\label{S5}
$
{\mathbf C}(A^\b \boxtimes B^\b)
\cong
\T(2)\otimes_{\T(1)\otimes \T(1)} \bigl({\mathbf C}(A^\b) \otimes {\mathbf
C}(B^\b)\bigr),
$
where $\T$ is the chain operad defined in Section \ref{s15}.
\end{theorem}

See \cite{MS4} for the proof.

The formula in Theorem \ref{S5} is precisely analogous to Definition V.1.1 of
\cite{KM}.  As a corollary to Theorem \ref{S5} we recover 
the results of \cite[Section V.3]{KM}, but with the
``linear isometries operad'' (actually the singular chains of the usual linear
isometries operad) replaced by $\T$.  The operad $\T$ has the advantage that it
is much smaller than  the linear isometries operad and its structure can be
described explicitly.











\begin{thebibliography}{99}

\bibitem{Bat93} Batanin, M.
Coherent categories with respect to monads and
coherent prohomotopy theory. Cahiers Topologie G\'eom. Diff\'erentielle
Cat\'egoriques 34 (1993),  279--304. 

\bibitem{Bat98} Batanin, M. 
Homotopy coherent category theory and $A\sb
\infty$-structures in monoidal categories. J. Pure Appl.\ Algebra 123 (1998),
67--103. 

\bibitem{Bat02} Batanin, M. 
The Eckmann-Hilton argument, higher operads and $E_n$-spaces.  Preprint 
available at http://front.math.ucdavis.edu/math.CT/0207281 

\bibitem{Bat03} Batanin, M.
The combinatorics of iterated loop spaces. Preprint available at 
http://front.math.ucdavis.edu/math.CT/0301221

\bibitem{BV1} 
Boardman, J.M. and Vogt, R.M. 
Homotopy-everything $H$-spaces. 
Bull.\ Amer.\ Math.\ Soc.\ 74 (1968), 1117--1122. 

\bibitem{BV2} 
Boardman, J.M.; Vogt, R.M. Homotopy invariant algebraic structures on
topological spaces. Lecture Notes in Mathematics, Vol.\ 347. Springer-Verlag, 
Berlin-New York, 1973.

\bibitem{CLM} Cohen, F.R., Lada, T.J., and May, J.P.
The homology of iterated loop spaces. 
Lecture Notes in Mathematics, Vol.\ 533. 
Springer-Verlag, Berlin-New York, 1976.

\bibitem{RCohen} Cohen, R.L. 
Stable proofs of stable splittings. 
Math.\ Proc.\ Cambridge Philos.\ Soc.\ 88 (1980), 149--151.

\bibitem{Curtis}
Curtis, E.B. Simplicial homotopy theory. Advances in Math.\ 6 (1971), 
107--209.

\bibitem{Deligne} Deligne, P. Letter to Stasheff et al. May 17, 1993.

\bibitem{EKMM}
Elmendorf, A.D., Kriz, I., Mandell, M.A., May, J.P. 
Rings, modules, and algebras in stable homotopy theory. 
With an appendix by M. Cole. Mathematical Surveys and Monographs, 47. 
American Mathematical Society, Providence, RI, 1997.

\bibitem{Gerst} Gerstenhaber, M.  The cohomology structure of an associative
ring.  Annals of Math.\, 78 (1963), 267--288.

\bibitem{Getzler} Getzler, E. Batalin-Vilkovisky algebras and two-dimensional
topological field theories.  Comm.\ Math.\ Phys.\ 159(1994), 265--285.

\bibitem{Grayson} Grayson, D.R.
Exterior power operations on higher $K$-theory.
$K$-Theory 3 (1989), 247--260.

\bibitem{HV} 
Hollender, J. and Vogt, R.M.
Modules of topological spaces, applications to homotopy limits and $E\sb
\infty$ structures. 
Arch.\ Math.\ (Basel) 59 (1992), 115--129.

\bibitem{HSS}
Hovey, M., Shipley, B. and Smith, J.H.
Symmetric spectra. 
J. Amer.\ Math.\ Soc.\ 13 (2000), 149--208.

\bibitem{HKV} Hu, P., Kriz, I. and Voronov, A. 
On Kontsevich's Hochschild Cohomology Conjecture. Preprint, September 2001. 

\bibitem{Kauf} Kaufmann, R. On spineless cacti, Deligne's conjecture and
Connes-Kreimer's Hopf algebra.  Preprint available at
http://front.math.ucdavis.edu/math.QA/0308005

\bibitem{K1} Kontsevich, M.
Operads and Motives in Deformation Quantization. Lett.\ Math.\ Phys.\  48
(1999) 35-72.

\bibitem{K2}
Kontsevich, M. and Soibelman, Y. Deformations of algebras over operads and the
Deligne conjecture. Conf\'erence Mosh\'e Flato 1999, Volume I, 255--307,
Math.\ Phys.\ Stud.\ 22, Kluwer Acad.\ Publ., Dordrecht, 2000.

\bibitem{KM} Kriz, I. and May, J.P.  Operads, algebras, modules and
motives.
Asterisque 233 (1995).


\bibitem{MacLane} Mac Lane, S. Categories for the working mathematician, 2nd 
Edition. Springer-Verlag, Berlin-New York, 1998.

\bibitem{MSSt} Markl, M., Shnider, S. and Stasheff, J. Operads in 
algebra, topology and physics.  Mathematical Surveys and Monographs, 96. 
American Mathematical Society, Providence, RI, 2002. 

\bibitem{MayG} May, J.P. The geometry of iterated loop spaces.  
Lectures Notes in Mathematics, Volume 271.
Springer-Verlag, Berlin-New York, 1972.

\bibitem{MayE} May, J.P. $E_\infty$ spaces, group completions, and 
permutative categories. New developments in topology (Proc.\ Sympos.\ 
Algebraic Topology, Oxford, 1972), pp.\ 61--93. London Math.\ Soc.\
Lecture Note Ser., No.\ 11, Cambridge Univ. Press, London, 1974. 

\bibitem{MS1} McClure, J.E. and Smith, J.H. A solution of Deligne's Hochschild
cohomology conjecture.  Recent progress in homotopy theory (Baltimore, MD,
2000), 153--193, Contemp.\ Math., 293, Amer.\ Math.\ Soc., Providence, RI,
2002.

\bibitem{MS2} McClure, J.E. and Smith, J.H.
Multivariable cochain operations and little n-cubes.  JAMS 16 (2003), 681-704.

\bibitem{MS3} McClure, J.E. and Smith, J.H. Cosimplicial Objects and little
n-cubes. I.  Preprint available at 
http://front.math.ucdavis.edu/math.QA/0211368

\bibitem{MS4} McClure, J.E. and Smith, J.H. 
Cosimplicial objects and little $n$-cubes. II. In preparation.

\bibitem{MS5} McClure, J.E. and Smith, J.H.  Cocylic objects and framed little
disks.  In preparation.

\bibitem{S} 
Sinha, Dev.
The topology of spaces of knots.  Preprint available at
http://front.math.ucdavis.edu/math.AT/0202287

\bibitem{Smith}
Smith, J.H. Simplicial group models for $\Omega\sp nS\sp nX$.  Israel 
J. Math.\ 66 (1989), 330--350.

\bibitem{Snaith} Snaith, V. P. 
A stable decomposition of $\Omega^n S^n X$. 
J. London Math.\ Soc.\ (2) 7 (1974), 577--583.

\bibitem{Stasheff} Stasheff, J.D.  Homotopy associativity of $H$-spaces. I.
Trans.\ Amer.\ Math.\ Soc.\ 108 (1963), 275--292.

\bibitem{T1} Tamarkin, D. Another proof of M. Kontsevich formality theorem.
Preprint available at http://front.math.ucdavis.edu/math.QA/9803025

\bibitem{T2} Tamarkin, D. Formality of Chain Operad of Small Squares.
Preprint available at http://front.math.ucdavis.edu/math.QA/9809164

\bibitem{V}
Voronov, Alexander A. Homotopy Gerstenhaber algebras. Conf\'erence
Mosh\'e Flato 1999, Vol.\ II (Dijon), 307--331, Math.\ Phys.\ Stud.\, 22,
Kluwer Acad.\ Publ., Dordrecht, 2000.

\bibitem{Weibel}
Weibel, C.A.
An introduction to homological algebra. 
Cambridge Studies in Advanced Mathematics, 38. 
Cambridge University Press, Cambridge, 1994. 

\end{thebibliography}
\end{document}